\documentclass[12pt]{article}

\usepackage{color}

\usepackage{amsmath}

\usepackage{graphicx}
\usepackage{float}
\usepackage{amsmath}
\usepackage{amsthm}
\usepackage{amssymb}
\usepackage{amscd}
\usepackage{amsfonts}
\usepackage{makeidx}
\usepackage{enumerate}      
\usepackage{IEEEtrantools}
\usepackage{mathrsfs}

\usepackage{amsmath,empheq}

\usepackage [all]{xy}
\usepackage{cases}

\usepackage{amsmath}

\usepackage{cancel}

\newtheorem{teo}{Theorem}[section]
    \newtheorem{lem}[teo]{Lemma}
    \newtheorem{prop}[teo]{Proposition}
    \newtheorem{coro}[teo]{Corollary}
    \newtheorem{defn}[teo]{Definition}
    \newtheorem{obs}[teo]{Remark}

\newcommand{\re}{\ref}

    \newtheorem*{dem}{\textsc{Proof}}
    \newcommand{\bdem}{\begin{dem}}
    \newcommand{\edem}{\end{dem}}
     \newcommand{\be}{\begin{equation}}
    \newcommand{\ee}{\end{equation}}
     \newcommand{\ba}{\begin{array}}
    \newcommand{\ea}{\end{array}}
\newcommand{\beqn}{\begin{eqnarray}}
    \newcommand{\eeqn}{\end{eqnarray}}
    \newcommand{\bl}{\begin{lem}}
    \newcommand{\el}{\end{lem}}
    \newcommand{\bp}{\begin{prop}}
    \newcommand{\ep}{\end{prop}}
\newcommand{\ds}{\displaystyle}

    \newcommand{\R}{\mathbb{R}}
    \newcommand{\C}{\mathbb{C}}
     \newcommand{\Z}{\mathbb{Z}}
     \newcommand{\no}{\noindent}

    \IEEEyessubnumber

\def\Re {{\rm Re\, }}                                       
 \def\Im {{\rm Im\,}}

\begin{document}
\title{Model BVP Problem for the Helmholtz equation in a nonconvex angle with periodic boundary data}
\author{
\large{ P.  Zhevandrov$^1$},\\
\large{A.  Merzon$^2$},\\
\large{M.I. Romero Rodr\'iguez$^3$}\\
\large{and J.E. De la Paz M\'endez$^4$}\\
%{\small{\it $^2$ Institute of Physics and  Mathematics, University of Michoac\'{a}n}},\\[-2mm]
%{\small 58060, Morelia, Michoac\'{a}n, M\'{e}xico}\\[-2mm]
%%%%%%%%%%%%%%%%%%%%%%%%%%%%%%%%%%%%%%%%%%%%%%%%%%%%%%%%%%%%%%%%
%%%%%%%%%%%%%%%%%%%%%%%%%%%%%%%%%%%%%%%%%%%%%%%%%%%%%%%%%%%%%%%%%
%%%%%%%%%%%%%%%%%%%%%%%%%%%%%%%%%%%%%%%%%%%%%%%%%%%%%%%%%%%%%%%%
%{\small}\\[-2mm]
{\small{\it $^1$ Facultad de  Ciencias F\'\i sico-Matem\'aticas, Universidad Michoac{a}na}},\\[-2mm]
{\small  Morelia, Michoac\'{a}n, M\'{e}xico}\\[-2mm]
{\small{\it $^2$ Instituto de F\'\i sica y  Matem\'aticas, Universidad Michoac{a}na}},\\[-2mm]
{\small  Morelia, Michoac\'{a}n, M\'{e}xico}\\[-2mm]
{\small{\it $^3$ Facultad de Matem\'aticas II, Universidad Aut\'onoma de Guerrero}},\\[-2mm]
{\small{Cd. Altamirano, Guerrero, M\'exico}}\\[-2mm]
{\small{\it $^3$ Facultad de Ciencias B\'asicas y Aplicadas, Universidad Militar Nueva Granada}},\\[-2mm]
{\small{Bogot\'a, Colombia}},\\[-2mm]
%%%%%%%%%%%%%%%%%%%%%%%%%%%%%%%%%%%%%%%%%%%%%%%%%%%%%%%%%%%%%%%%%%%%%
{\small{\it E-mails}: pzhevand@gmail.com,}
{\small anatolimx@gmail.com,}\\[-2mm] {\small jeligio12@gmail.com,}{\small maria.romeror@unimilitar.edu.co} }
\maketitle

\begin{abstract}
 In the presented  work, we solve the Dirichlet boundary  problem  for the Helmholtz equation in an exterior  angle with periodic boundary data. We prove the existence and uniqueness of solution in an appropriate funcional class and we give an explicit  formula for it  in the form of the Sommerfeld  integral. The method of complex characteristics \cite {BKM} is used.
\end{abstract}

\section{Introduction}
\setcounter{equation}{0}

\begin{figure}[htbp]
\centering
\includegraphics[scale=0.27]{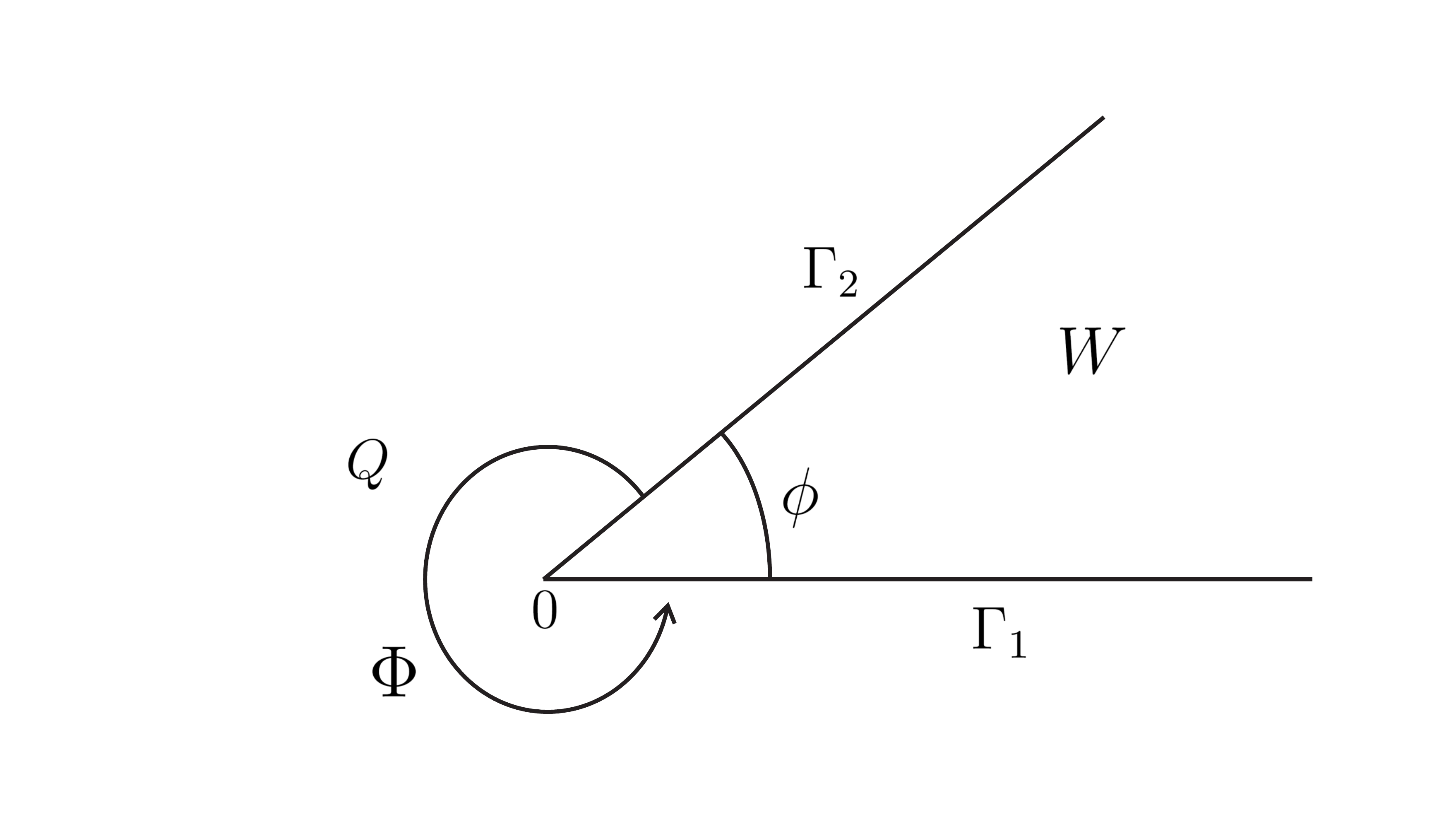}
\caption{BVP in exterior angle}\label{AQ}
\end{figure}

We consider the following model boundary value problem (BVP) for the Helmholtz equation in a plane angle $Q$ of magnitude $\Phi>\pi$ with a complex frequency $\omega\in\C^{+}=\big\lbrace \Im\omega>0\big\rbrace$:
\begin{equation}\label{D}
\left\{     \begin{array}{rcl}
    (-\Delta-\omega^2)u(x) =0, &  & x\in Q \\\\
              B_1 u(x)\Big\vert_{\Gamma_1} = f_1(x), \\\\
              B_2 u(x)\Big\vert_{\Gamma_2}=f_2(x).   &  &  
            \end{array} 
\right.
\end{equation}
Here $\Gamma_{l}$ for $l=1,2$  are the sides of the angle $Q$, $B_{l}=I$ or $B_{l}=\ds\frac{\partial}{\partial n_{l}}$ ($n_l$ is the exterior normal to $\Gamma_l$), $f_{l}$ are given functions which can be distributions, see Fig.\;\ref{AQ} 

\newpage

Problems of this type arise is many areas of mathematical physics.
 We list some of them.\\ 
 Firstly, such BVP describe the diffusion of a desintegrating gas \cite[Ch. VII, cf. 2]{Tikh}.\\
Secondly, diffraction problems by the wedge $W=\R^2\setminus Q$ are reduced to such problems for a lossy medium \cite{Teix 1991} or for a slightly conducting medium \cite{ST 1989}.\\
Thirdly, time-dependent diffraction by wedges \cite{Sob34,KB,Ka,Somm,obe,BER JML1,Ro,Rot,kmjv}
 is reduced to a problem of this type after the Fourier-Laplace (F-L) transform $t\to \omega$ with respect to time \cite{KMM, km, mm, mesqNN}.\\
Fourthly, the problem of scattering of waves emitted by a point source by $W$ is reduced
to this problem, with  
$$
f_{1}(y_1)= e^{ik y_1},\quad y_1>0,\quad f_2=0.
$$
Let us describe this scattering problem in more detail  since it was precisely this problem which was the starting point of the present paper.\\

\begin{figure}[htbp]
\centering
\includegraphics[scale=0.3]{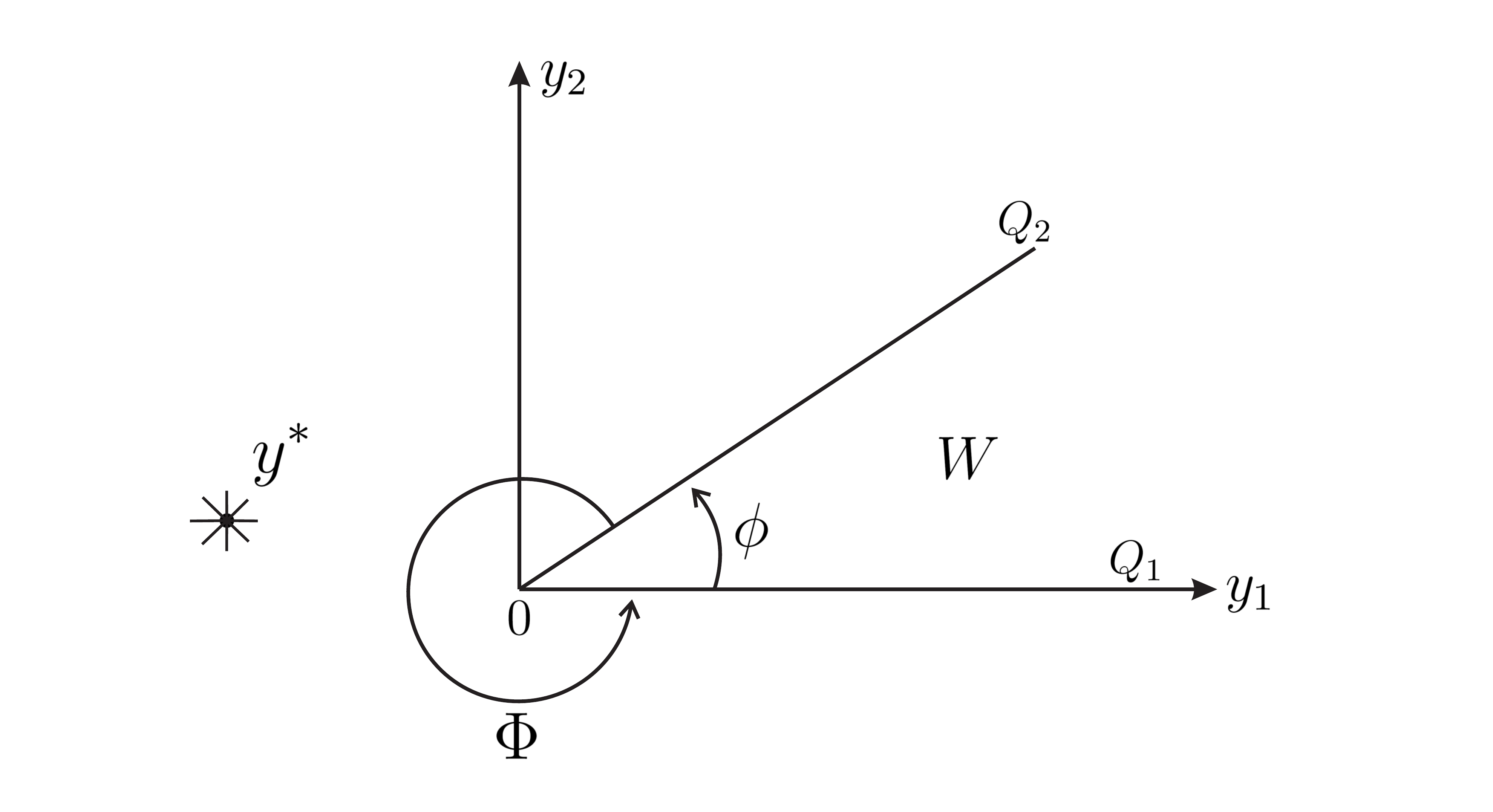}
\caption{}\label{W}
\end{figure}

\medskip
\no Let us consider the following time-dependent scattering problem:
\begin{equation}\label{we}
\left.
\begin{array}{rcl}
U_{tt}(y,t) =\Delta U(y,t)+\delta(y-y^{\ast}) e^{-i\omega_0 t},\quad y\in Q:=\R^2\setminus W,\quad t\geq 0,\quad \omega_0\geq0\\\\
U(y,t)\Big\vert_{\partial Q} =0.
\end{array}
\right|
\end{equation}
Here $y^{\ast}\in Q$.
After the F-L transform $t\to\omega: U(y,t)\longrightarrow \hat{U}(y,\omega)$ (\ref{we}) becomes equivalent to
\begin{equation}\label{fwe}
\left.
\begin{array}{rcl}
-\omega^2\hat{U}(y,\omega)&=&\Delta\hat{U}(y,\omega)+\delta(y-y^{\ast})\ds\frac{i}{\omega-\omega_0},\quad\omega_0>0\\\\
\hat{U}(y,\omega)\Big\vert_{\partial Q}&=&0,
\end{array}
\right| \omega\in\C^{+}
\end{equation}
see Fig.\;\ref{W}

\newpage

We reduce this problem to a homogeneous Helmholtz equation with nonhomogeneous boundary conditions.
Let $\mathscr{E}(y,\omega)$ be s.t.
\begin{equation}\label{ge}
-\Delta\mathscr{E}(y,\omega)-\omega^2\mathscr{E}(y,\omega)=\frac{i\delta(y-y^{\ast})}{\omega-\omega_0},\quad y\in\R^2,\quad \Im\omega>0,\quad\omega_0>0.
\end{equation}
Passing to the Fourier transform $y\to\eta$ in (\ref{ge}), we obtain
\begin{equation*}
\tilde{\mathscr{E}}(\eta,\omega)=\frac{i e^{i\eta\cdot y^{\ast}}}{(|\eta|^2-\omega^2)(\omega-\omega_0)}.
\end{equation*}
Hence
\begin{equation*}
\mathscr{E}(y,\omega)=\frac{i}{(\omega-\omega_0)(2\pi)^2}\int\limits_{\R^2}\frac{e^{-i\eta(y-y^{\ast})}}{|\eta|^2-\omega^2}\;d\eta.
\end{equation*}
Define
\begin{equation*}
U_{s}(y,\omega):=\hat{U}(y,\omega)-\mathscr{E}(y,\omega),\quad y\in Q,\quad \omega\in\C^{+}.
\end{equation*}
By (\ref{fwe}), (\ref{ge}) this implies that $U_s(y,\omega)$ satisfies the problem of type (\ref{D}):
\begin{equation}\label{su_s}
\left\{     \begin{array}{rcl}
              -\Delta U_{s}(y,\omega)-\omega^2 U_{s}(y,\omega) = 0, &  & y\in Q \\\\
              U_{s}(y,\omega)\Big\vert_{\partial Q} = -\mathscr{E}(y,\omega)\Big\vert_{\partial Q} &  & 
            \end{array} 
\right|    \ \omega\in\C^{+}.
\end{equation} 
Let us calculate $\mathscr{E}(y,\omega)$ on $\partial Q$. We have 
\begin{equation*}
-\mathscr{E}(y, \omega)\Big\vert_{Q_1}=-\mathscr{E}(y_1, 0, \omega)= C(\omega)\ds\int\limits_{\R^2} e^{-i\eta_1 y_1}\; C_1(\eta, \omega,y^{\ast})\;d\eta
\end{equation*} 
where $C(\omega)=\ds\frac{i}{(\omega-\omega_0)(2\pi)^2}$; $C_1(\eta, \omega, y^{\ast})= \ds\frac{e^{i\eta\cdot y^{\ast}}}{|\eta|^2-\omega^2}$. 
Hence,\\
 $-\mathscr{E}(y, \omega)\Big\vert_{Q_2}= C(\omega)\ds\int\limits_{\R^2} e^{-i(\eta_1\cot\Phi+\eta_2)y_2}\; C_1(\eta, \omega,y^{\ast})\;d\eta$.
Thus (\ref{su_s}) is equivalent to
\begin{equation}\label{pU_s}
\left\{     \begin{array}{rcl}
              & &-\Delta U_{s}(y,\omega)-\omega^2 U_{s}(y,\omega) = 0, \quad y\in Q \\\\
              & & U_{s}(y_1,0,\omega) = C(\omega)\ds\int\limits_{\R^2} e^{-i\eta_1 y_1}\;C_1(\eta, \omega, y^{\ast})\;d\eta, \quad y_1>0\\\\
& & U_{s}(y_2\cot\phi, y_2,\omega)= C(\omega) \ds\int e^{-i(\eta_1\cot\Phi+\eta_2)y_2}\;C_1(\eta, \omega, y^{\ast})\;d\eta, \quad y_2>0.               
            \end{array} 
\right. 
\end{equation}
Therefore it seems natural to solve first the following model problem, corresponding to problem (\ref{pU_s})
\begin{equation}\label{PU}
\left\{     \begin{array}{rcl}
              & &-\Delta U(y,\omega)-\omega^2 U(y,\omega) = 0, \quad y\in Q \\\\
              & & U(y_1,0) = e^{-ik_1 y_1}, \quad y_1>0\\\\
& & U(y_2\cot\phi, y_2)= e^{\frac{-ik_2 y_2}{\sin\Phi}}, \quad y_2>0,               
            \end{array} 
\right. 
\end{equation}
where $k_1, k_2\in\R$.\\ 
Note that this problem is similar to the BVP in \cite[(23)]{KMM} arising in the problem of time-dependent diffraction:
\begin{equation}\label{23}
\left\{     \begin{array}{rcl}
              & &\big(-\Delta-\omega^2\big) U(\omega,y) = 0, \quad y\in Q \\\\
              & & U(\omega, y) = -g(\omega) e^{i\omega y_1\cos\alpha}, \quad y\in Q_1\\\\
& & U(\omega, y)= -g(\omega) e^{-i\omega y_2 \big(\cos(\alpha+\Phi)/\sin\Phi\big)}, \quad y\in Q_2.               
            \end{array} 
\right. 
\end{equation}
However, there are substantial  differences.  The exponents in the right-hand side of
this problem are complex ($\Im \omega>0$) and, thus, the corresponding functions decrease exponentially when $y_{1,2}\longrightarrow +\infty$ since $\alpha<\phi<\ds\frac{\pi}{2}$.\\ Moreover the structure of boundary conditions in (\ref{23}) is connected with the first equation through the common parameter $\omega$. This gives a unique opportunity to reduce problem (\ref{23}) to a difference equation which is solved easily in an explicit form.\\
In contrast, problem (\ref{PU}) has periodic boundary conditions which are independent of the first equation. This results in the fact that the corresponding difference equation cannot be solved as easily as in the previous case except for the case when $\Phi=\ds\frac{3\pi}{2}$ (see Section\;5).\\
In turn, problem (\ref{PU}) splits  into two problems for $u_1, u_2$ such that $U=u_1+u_2$ by linearity:
\begin{equation}\label{u_1}
\left\{     \begin{array}{rcl}
              & &-\Delta u_1(y,\omega)-\omega^2 u_1(y,\omega) = 0, \quad y\in Q \\\\
              & & u_1(y_1,0) = e^{-ik_1 y_1}, \quad y_1>0\\\\
& & u_1(y_2\cot\phi, y_2)= 0, \quad y_2>0.               
            \end{array} 
\right. 
\end{equation}
\medskip

\begin{equation}\label{u_2}
\left\{     \begin{array}{rcl}
              & &-\Delta u_2(y,\omega)-\omega^2 u_2(y,\omega) = 0, \quad y\in Q \\\\
              & & u_2(y_1,0) = 0, \quad y_1>0\\\\
& & u_2(y_2\cot\phi, y_2)= e^{-\frac{ik_2 y_2}{\sin\phi}}, \quad y_2>0.               
            \end{array} 
\right. 
\end{equation}
Here $\omega\in\C^{+}$.
This paper is devoted to solving the model problem (\ref{u_1}). Solution of (\ref{u_2}) is obtained from (\ref{u_1}) by a simple change of variable, (see (\ref{u_1''})).\\
Note that the BVP in a right angle $Q$ or in its complement and in other particular angles whose magnitudes are commensurate with $\pi$, were considered in many papers \cite{Mei87, MPR,  Mei93,  Mei94, ms, Pen99, Cas08, CK08, caskap, 4, Cas06}.\\
In those papers exact results were obtained by means of operator methods. Boundary data in those papers belong to Sobolev spaces $H_s(\R)$, $s>0$. We consider another type of boundary data, namely, periodic functions. We obtain exact solutions in explicit form, namely, in the form of Sommerfeld type integrals. We use the method of automorphic functions (MAF) on complex characteristics  \cite{BKM}. This method was developed by A.Komech  for $\Phi<\pi$ in \cite{K} and then was extended to $\Phi>\pi$ in \cite[Section\;1.2 and part\;2]{BKM}. It allows us to find all distributional solutions of the BVP for the Helmholtz equation in arbitrary angles with general boundary conditions. It was applied, in particular, to time-dependent diffraction problems by angles \cite{mm, KMJ2014, MZJ2015, KM14, KM15, KMNMV2018}.\\      
It should be noted that there is a very effective Sommerfeld-Malyzhinetz method of constructing solutions of diffraction  problems in angles; by means of this method many important results were obtained \cite{BLG}. This method allows one to obtain the solution in the form of the Sommerfeld integral. However, this method does not allow one to prove uniqueness which usually is proved on the basis of physical considerations.

We also obtain solution in the Sommerfeld integral form, using the MAF which also allows us to prove uniqueness in an appropriate functional space (see e.g. \cite{Mei87}).

The paper is organized as follows: in Section 2 we formulate the main result, the Sections 3-10 are devoted to its proof. In Section 3 we reduce boundary value problem to a difference equation and we prove the necessary and sufficient  conditions for the existence of solution. In Sections 4 and 5 we find  solution of the difference equation for $\Phi\not = \frac{3}{2}\pi$ and $\Phi\not = \frac{3}{2}\pi$  respectively. In Section 6 we prove the asymptotics of the integrand for the Sommerfeld's representation of the solution. In Sections 7 and 8 we give the Sommerfeld-type representation of solution and we prove the boundary conditions. In Section 9 we prove the existence and uniqueness of the solution. En Appendices we prove some technical assertions.

\section{The main result}\label{MR}
\setcounter{equation}{0}

We will construct the solutions of problems (\ref{u_1}) and (\ref{u_2}) in the form of the well-known Sommerfeld integrals which have the form 
\begin{equation*}
\ds\int\limits_{\mathcal{C}} e^{-\omega\rho\sinh w} v(w+i\theta)\;dw, 
\end{equation*}
where $\mathcal{C}$ is a certain contour on the complex plane and the correct construction of the factor $v(w)$, which ensures that (\ref{SI}) satisfies the boundary conditions, is the main difficulty of the problem.\\ 
To formulate the main result we need to describe the integrand $v(w)$ of the Sommerfeld integral. The construction of this integrand is the main contents of this paper.\\
Consider $\hat{v}_1(w)$ given by (\ref{check{v}_1}) where $\hat{v}_1^1(w)$ is given by (\ref{hat{v}_1^1}) for  $\Phi\neq \ds\frac{3\pi}{2}$ and by (\ref{u_3}) for $\Phi=\ds\frac{3\pi}{2}$ with $\hat{G}$, given by (\ref{G(w)}).
Let $(\rho,\theta)$ be the polar coordinates in $\R_{y}^{2}$,
\begin{equation}\label{PC}
y_1=\rho\cos\theta,\quad y_2=\rho\sin\theta, \quad \rho>0,\quad \theta\in\big[2\pi-\Phi,2\pi\big],              
\end{equation}
and $\mathcal{C}$ be the Sommerfeld double-loop contour, see Fig.\;\ref{CS3}.  
\begin{defn}
$E$ is the space of the $C^{\infty}\big(\overline{Q}\setminus\lbrace 0\rbrace\big)$ functions bounded 
with its first derivatives  in $\overline{Q}\setminus B_{\varepsilon}(0)~\forall \varepsilon>0$ and  admitting the following asymptotics at the origin
\begin{equation}\label{SP}
\left.  \begin{array}{rcl}
& & u(\rho,\theta)= C(\theta)+o(1)\\\\
& & \nabla u(\rho,\theta)= C_1(\theta) \rho^{-1}+C_2(\theta)+o(1)
\end{array}
\right|  \quad \rho\to0.
\end{equation}
\end{defn}
Our main result consists in the following statement.
\begin{teo}\label{Tu_1u_2}
{\bf i)} Let $\omega\in\C^{+}, k_1, k_2>0$. There exists a solution to problem (\ref{PU})
\begin{equation*}
U(\rho, \theta)=u_1(\rho, \theta)+u_2(\rho, \theta),\quad (\rho, \theta)\in Q,
\end{equation*}
belonging to $E$, where $u_1, u_2$ are solutions to (\ref{u_1}), (\ref{u_2}), respectively, which admit a Sommerfeld integral representation   
\begin{equation}\label{u_1'}
u_1(\rho,\theta)=\ds\frac{1}{4\pi\sin\Phi}\int\limits_{\mathcal{C}} e^{-\omega\rho\sinh w}\hat{v}_1(w+i\theta) dw, 
\end{equation}
\begin{equation}\label{u_1''}
u_2(\rho,\theta)=\ds\frac{1}{4\pi\sin\Phi}\int\limits_{\mathcal{C}} e^{-\omega\rho\sinh w}\hat{v}_1(w+i\theta_1) dw, 
\end{equation}
where in (\ref{u_1'}) $\hat{v}_1$ is constructed according to the algorithm presented below for $k=k_1$, in (\ref{u_1''}) for $k=k_2$  and 
\begin{equation}\label{Ch}
\theta_1= -\theta+4\pi-\Phi.
\end{equation}
{\bf ii)} The solution $U$ is unique in $E$. 

\end{teo}
\begin{obs}\label{rbc}
The integrals in (\ref{u_1'}), (\ref{u_1''}) converge absolutely since the infinite part of $\mathcal{C}$ belongs to the region of the superexponential decrease of  $e^{-\omega \rho\sinh w}$ (see Fig.\;\ref{C_1C_2}) and $\hat{v}_1$ admits asymptotics (\ref{as hat{v}_1}). The boundary values are understood in the sense of distributions. 
\end{obs}
\begin{obs}\label{r22}
The representation (\ref{u_1''}) follows easily from (\ref{u_1'}) by the change (\ref{Ch}).
\end{obs}
\begin{obs}\label{ast1}
Note that the solution $u_{1}$ also admits a slightly different representation where several different Sommerfeld-type contours are used (see (\ref{des})). 
\end{obs}

\section{Reduction to a difference equation. Necessary and sufficient conditions for the Neumann data}
\setcounter{equation}{0}

Consider problem (\ref{u_1}). The MAF permits to reduce this problem to finding the Neumann data of solution $u_1$, and it consists of several steps. In the following subsections we present these steps.\\ 
We assume that the solution $u_1\in S'(Q):=\Big\lbrace u\Big\vert_{Q}, u\in S'(\R^2)\Big\rbrace$.
The first step of the MAF  is to reduce the problem to the complement of the first quadrant  and to extend the solution $u_1$ to the plane, see \cite{BKM, km}.
 
\subsection{First step: extension of $v_{l}^{\beta}(x_{l})$ to the whole plane $\R^2$}

Consider the linear transformation
\begin{equation*}
\mathcal{J}(y): x_1=y_1+y_2\cot\Phi,\quad x_2=-\frac{y_2}{\sin\Phi}, 
\end{equation*}
which sends the angle $Q$ to the right angle $K:=\Big\lbrace (x_1,x_2): x_1<0~{\rm or}~x_2<0\Big\rbrace$. This transformation reduces system (\ref{u_1}) to the  problem (\ref{ra1})-(\ref{ra3}) in the complement $K$ of the first quadrant for
\begin{equation*}
v(x):=u_1\big(\mathcal{J}^{-1}(y)\big)
\end{equation*}

\begin{subequations}
\begin{empheq}[left=\empheqlbrace]{align}
\mathscr{H}(D) v(x) &= 0, \qquad\quad~ x\in K 
\label{ra1}
\\ \nonumber\\
v(x_1,0) &= e^{-ik x_1}, \quad x_1>0
\label{ra2}
\\ \nonumber\\
v(0,x_2) &= 0,   \qquad\quad ~x_2>0,
\label{ra3} 
\end{empheq}
\end{subequations}
where 
\begin{equation}\label{mathcal}
\mathscr{H}(D)=-\ds\frac{1}{\sin^2\Phi}\Bigg[\Delta-2\cos\Phi\ds\frac{\partial^2}{\partial x_1\partial x_2}\Bigg]-\omega^2.
\end{equation}
By \cite[Lemma\;8.2]{BKM}, if $v(x)\in S'(K)$ is a solution of equation (\ref{ra1}), then there exists an extension $v_0\in S'(\R^2)$ of $v$ by 0,  such that $v_0\Big\vert_{K}=v$,
 \begin{equation}\label{Pu_0}
\mathscr{H}(D)v_0(x)=\gamma(x),\quad x\in\R^2,
\end{equation}
where $\gamma\in S'(\R^2)$ and has the form
\begin{equation}\label{gamma_a}
\begin{array}{lll}
\gamma (x)&=&\ds\frac{1}{\sin^2\Phi}\Big[\delta(x_2)v_1^1(x_1)+\delta'(x_2)v_1^0(x_1)+\delta(x_1)v_2^1(x_2)+\delta'(x_1)v_2^0(x_2)-\\\\
&-& 2\cos\Phi\;\delta(x_2)\;\partial_{x_{1}} v_1^{0}(x_1)-2\cos\Phi\;\delta(x_1)\partial_{x_2}\;v_2^{0}(x_2)\Big],\quad x\in\R^2;
\end{array}
\end{equation}
$v_{l}^{\beta}(x_{l})\in S'(\overline{\R^{+}}):=\Big\lbrace v\in S'(\R): {\rm supp}\;v\subset \overline{\R^+}\Big\rbrace$.\\
We will use the extension  of the Fourier transform $F$ defined on $S(\R)\subset S'(\R^2)$,  $ \varphi(x_1,x_2)\to \tilde{\varphi}(z_1,z_2)$, $\varphi\in S(\R^2)$ to $S'(\R^2)$ by continuity:
\begin{equation}\label{ft}
F_{x\to z}\big[\varphi\big](z)=F\big[\varphi(x)\big](z)=\tilde{\varphi}(z_1,z_2):=\iint e^{iz_1x_1+iz_2 x_2} \varphi(x_1, x_2)\;dx_1 dx_2,
\end{equation}
and denote this extension by the tilde,  $\tilde{v}(z)=F_{x\to z}[v(x)],~v\in S'(\R^2)$. Applying this transform to (\ref{Pu_0}) and using the fact that  $\mathscr{H}(z)\neq 0, z\in\R^2$, we obtain 
\begin{align}\label{div1}
\tilde{v}_0(z)&=\ds\frac{\tilde{\gamma}(z)}{\mathscr{H}(z)},\quad z\in\R^2,\quad \tilde{v}_0\in S'(\R^2). 
\end{align}
Hence,
\begin{equation}\label{8.5}
v_0(x)=F^{-1}_{x\to z}\Bigg[\ds\frac{\tilde{\gamma}(z)}{\mathscr{H}(z)}\Bigg],\quad x\in\R^2,\quad v(x)=v_0(x)\Big\vert_{K}.
\end{equation}
Here $\tilde{\gamma}(z)$ is the Fourier transform of (\ref{gamma_a}), and for $z\in\R^2$
\begin{equation}\label{tilde{gamma}}
\begin{array}{lll}
\tilde{\gamma} (z)
= \ds\frac{1}{\sin^2\Phi}\Big[\tilde{v}_1^1(z_1)-\tilde{v}_1^0(z_1)\big(iz_2-2\cos\Phi\;iz_1\big)+\tilde{v}_2^1(z_2)-
 \tilde{v}_2^0(z_2)\big(iz_1-2\cos\Phi\;iz_2\big)\Big],
\end{array}
\end{equation}
where $\tilde{v}_{l}^{\beta}(z_l)$ are the Fourier transforms of $v_{l}^{\beta}(x_{l})$. Thus, if we know $v_{l}^{\beta}(x_{l})$, we know $v$ by (\ref{8.5}), and   problem (\ref{gamma_l^beta}) is reduced to finding the four  functions $\tilde{v}_{l}^{\beta}(x_{l})$, $l=1,2$, $\beta=0,1$.
\begin{obs}\label{rl}
Formula (\ref{gamma_a}) is obtained by direct  differentiation (in the sense of 
$\mathcal{D}'(\R^2))$ of the discontinuous function
\begin{equation}\label{ex}
v_0(x)=\Big[v(x)\Big]_0:=\left\{
\begin{array}{rcl}
v(x), &  & x\in K\\\\
0,     &  & x\notin K
\end{array}
\right.
\end{equation}
in the case when $v(x)\in C^{\infty}(\overline{K})$. Moreover, the formula
$$
f'_{0}(x)=\Big[f'(x)\Big]_{0}+f(0)\delta(x),\quad x\in\R
$$
is used for $f\in C^{\infty}\Big(\overline{\R^{+}}\Big)$. Obviously, in this case the functions $v_l^{\beta}$ in (\ref{gamma_a}) are the Cauchy data of the function $v$:
\begin{equation}\label{vlb}
\begin{array}{lll}
v_1^0(x_1)= v(x_1,0),~ x_1>0,  \qquad \quad  v_1^1(x_1) = \ds\frac{\partial}{\partial x_2} v(x_1,0-),\quad~  x_1>0\\\\
 v_2^0(x_2) = v(0-,x_2), ~ x_2>0,\qquad    v_2^1(x_2) =\ds\frac{\partial}{\partial x_1} v(0-,x_2), ~\quad x_2>0.
\end{array}
\end{equation}
It turns out that formula (\ref{gamma_a}) and representations (\ref{vlb}) remain true for distributional solutions. The following two lemmas describe the solution of equation (\ref{ra1}) in terms of its Cauchy data.
\end{obs}

\begin{lem}\label{8.3}
\cite[Lemma\;8.3]{BKM}. Let $v\in S'(K)$ be a distributional solution to equation (\ref{ra1}) and let $v_0$ be its extension by 0 to $\R^2$ satisfying (\ref{Pu_0}), (\ref{gamma_a}). Then the Cauchy data
\begin{equation*}
\left\{
\begin{array}{rcl}
              v_1^{\beta}(x_1):=\partial_2^{\beta} v_0(x_1,0-) &  & x_1>0 \\\\
              v_2^{\beta}(x_2):=\partial_2^{\beta} v_0(0-,x_2) &  & x_2>0
            \end{array} 
\right|    \ \beta=0,1
\end{equation*}
exist (here the limits are understood in the sense  $\mathcal{D}'(\R^+):=\mathcal{D}'(\R)\Big\vert_{\R^{+}}).$
\end{lem}
\begin{lem}\label{l8.4}
\cite[Lemma\;8.4]{BKM}. Let $v\in S'(K)$ be a distributional solution to equation (\ref{ra1}) given by (\ref{8.5}), where $\gamma$ is defined  by (\ref{gamma_a}). Then  $v_l^{\beta}\Big\vert_{\R^{+}}$ are the Cauchy data of $v$. 
\end{lem}
\begin{obs}
Formula (\ref{8.5}) and Lemma\;\ref{l8.4} show that it suffices to find the Neumann data $v_1^1, v_2^1$ in (\ref{gamma_a}) in order to solve problem (\ref{ra1})-(\ref{ra3}).  
\end{obs}

\medskip

Now we use  boundary conditons (\ref{ra2}), (\ref{ra3}). Let $v\in S'(K)$ be a solution to (\ref{ra1})-(\ref{ra3}) and $v_0\in\mathcal{S}'(\overline{K})$ be its extension by $0$ satisfying (\ref{Pu_0}); then,
 by (\ref{vlb}).
\begin{equation}\label{gamma_l^beta}
\left.
\begin{array}{rcl}
v_1^0(x_1)\Big\vert_{\R^+}= e^{-i k x_1},\quad x_1>0,\qquad
v_2^0(x_2)\Big\vert_{\R^+}= 0,\quad x_2>0.
\end{array}
\right.
\end{equation}
Since ${\rm supp}~v_{l}^{0}(x_{l})\subset \overline{\R^+}$, by the distribution theory, we have, generally speaking, 
\begin{equation}\label{deltam}
\left\{
\begin{array}{rcl}
v_1^0(x_1)&=&\Big[e^{-ik x_1}\Big]_0+c_1^0\delta(x_1)+c_1^1\delta'(x_1)+\cdots +c_1^{m}\delta^{(m)}(x_1),\\\\
v_2^0(x_2)&=&c_2^0\delta(x_2)+c_2^1\delta'(x_2)+\cdots +c_2^{m}\delta^{(m)}(x_2),
\end{array}
\right.
\end{equation}
for some $m\geq0$.
Here $\Big[e^{-ikx_1}\Big]_0$ is defined similarly to (\ref{ex}). Obviously $\Big[e^{-ikx_1}\Big]_0=\Theta(x) e^{-ikx_1}$ where $\Theta(x)$ is the Heaviside function.
We will find a solution to (\ref{ra1})-(\ref{ra3}) for  
\begin{equation}\label{as}
c_1^0=\cdots =c_1^m=c_2^0=\cdots=c_2^m=0.
\end{equation}
Thus we put 
\begin{equation}\label{v_1^0}
v_1^0(x_1)= \Big[e^{-ik x_1}\Big]_0=\left\{
\begin{array}{rcl}
e^{-ikx_1},& & x_1\geq0\\\\
0,& &  x_1<0
\end{array}
\right. \ \in S'(\R),\quad v_2^0(x_2)=0.
\end{equation}

\begin{obs}
It is not guaranteed a priori that the solution of (\ref{ra1})-(\ref{ra3}) exists under condition (\ref{as}) because $\tilde{v}_{l}^{\beta}$ should satisfy a certain connection equation (see Section\;3.2). Nevertheless, it turns out that we are able to construct an explicit solution under the condition (\ref{as}). Solutions which correspond to nonzero values of $c_{l}^m$ in (\ref{deltam}) will only contain additional singularities at the origin, and are not of interest.
\end{obs}
\no Substituting $v_1^0, v_2^0$ given by (\ref{v_1^0}) in (\ref{Pu_0}), we obtain
\begin{equation}\label{e3}
\mathscr{H}(D) v_0(x)=\gamma(x),
\end{equation}
with $\gamma$ containing only two unknown functions $v_1^1$ and $v_2^1$.\\ 
The MAF gives the necessary and sufficient conditions for the functions $v_1^1$
and $v_2^1$, which allow us to find these functions in an explicit form. Substituting these functions in (\ref{e3}) we obtain $v_0$ (and so $v$) by (\ref{8.5}), (\ref{gamma_a}).\\
In what follows we consider equation (\ref{e3}).

\subsection{Second step: Fourier-Laplace transform and the lifting to the Riemann surface. Connection equation}

\begin{figure}[htbp]
\centering
\includegraphics[scale=0.32]{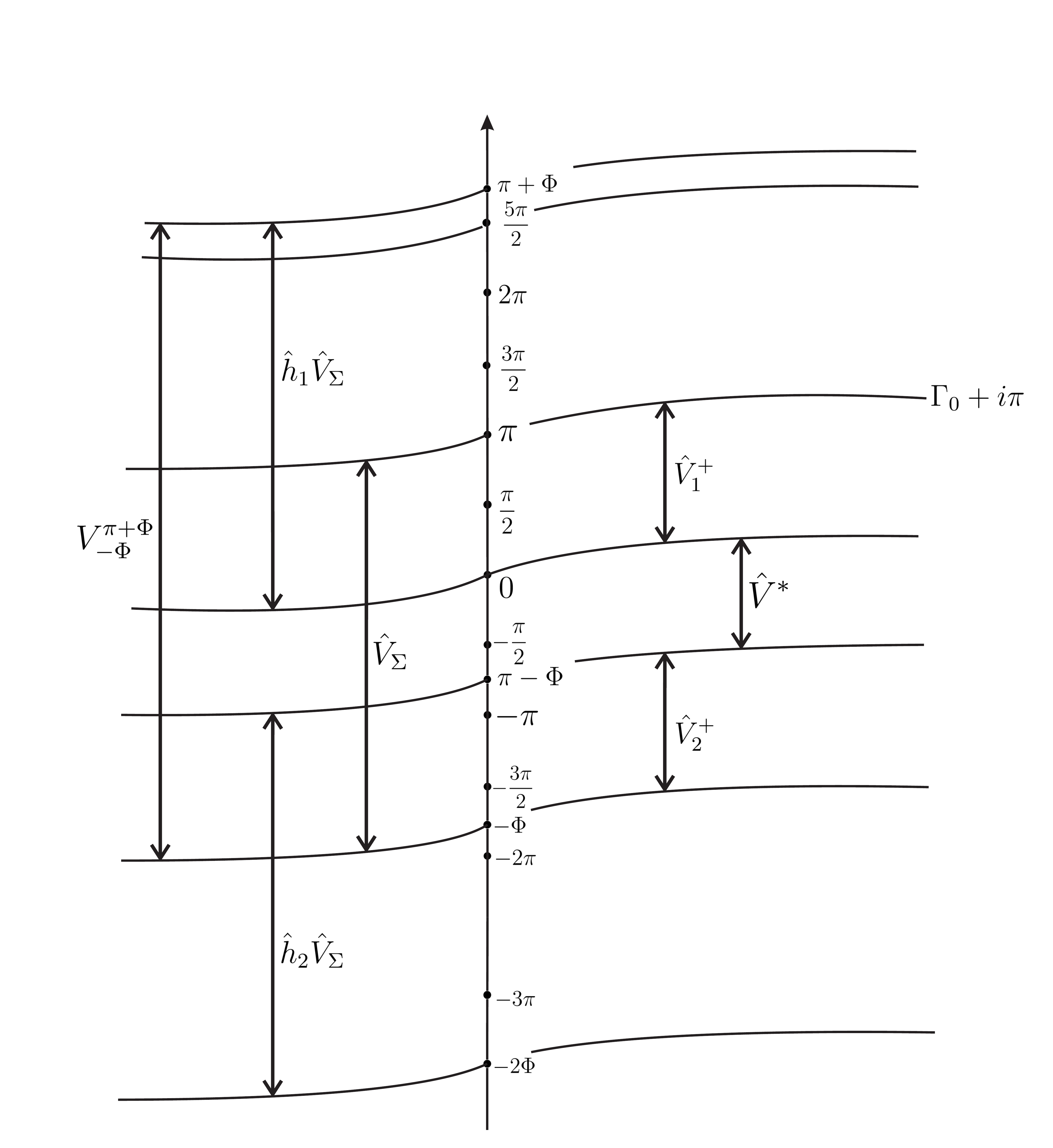}
\caption{Connection Equation}\label{w=i}
\end{figure}

In addition to the (real) Fourier transform (\ref{ft}) we will use the complex Fourier transform (or Fourier-Laplace (F-L) transform). Let 
\begin{equation*}
f\in S'(\overline{\R^{+}}):=\Big\lbrace f\in S'(\R): {\rm supp}\;f\subset \overline{\R^{+}}\Big\rbrace.
\end{equation*}
Then by the Paley-Wiener theorem \cite{RS},  $\Big({\rm see}~{\rm also}~\cite[Theorem\;5.2]{Ka}\Big)$ $\tilde{f}(z)=F\big[f\big]\in\R$, admits an analytic continuation $\tilde{f}(z)\in\mathcal{H}\big(\C^+\big),~ \C^{+}:=\Big\lbrace z\in\C: \Im z>0\Big\rbrace$ and $\lim \tilde{f}(z_1+iz_2)=\tilde{f}(z_1)~{\rm in}~S'(\R),~ {\rm as}~\varepsilon\to 0+$. Since $v_{l}^{\beta}(x_l)\in S'(\overline{\R^{+}})$, there exist their F-L transforms 
\begin{equation}\label{atilde{v}_l}
\tilde{v}_{l}^{\beta}(z_l)\in\mathcal{H}\big(\C^+\big),\quad l=1,2;\quad \beta=0,1. 
\end{equation}
In particular,
from (\ref{v_1^0}) we have
\begin{equation}\label{tilde{v}_1^0}
\tilde{v}_1^0(z_1)=\frac{i}{z_1-k},\quad z_1\in\overline{\C^{+}}, 
\end{equation}
where for $z_1\in\R$, $\tilde{v}_1^{0}(z_1)=\ds\lim_{\tau_1\to 0+}\tilde{v}_1^0(z_1+i\tau_1)$ in $S'(\R)$.
Hence, using (\ref{tilde{gamma}}) we obtain (since  $v_2^0\equiv 0$) 
\begin{equation}\label{tilde{gamma}''}
\tilde{\gamma}(z)=\ds\frac{1}{\sin^2\Phi}\Bigg[\tilde{v}_1^1(z_1)+\ds\frac{z_2-2\cos\Phi\;z_1}{z_1-k}+\tilde{v}_2^1(z_2)\Bigg],\quad z\in\R^2.
\end{equation}
In the MAF, the Riemann surface of complex zeros of the symbol of the operator (\ref{mathcal})  plays an essential role, since a necessary condition for  the existence  of the solution on $\tilde{\gamma}(z)$ can be written in terms of this surface. The symbol of this operator is the polynomial 
$$
\mathscr{H}(z)=\ds\frac{1}{\sin^2\Phi}\Big(z_1^2+z_2^2-2z_1 z_2\cos\Phi\Big)-\omega^2,\quad (z_1,z_2)\in\C^2.
$$
Obviously, $\mathscr{H}(z)$ does not have real zeros, but it does have complex ones. Denote the Riemann surface of the complex zeros of $\mathscr{H}$ by
\begin{equation*}
V:=\Big\lbrace z\in\C^2: \mathscr{H}(z)=0\Big\rbrace.
\end{equation*} 
It is convenient to parametrize the complex surface $V$ introducing the parameter $w\in\C$.\\
The Riemann surface $V$ admits a universal covering $\hat{V}$, which is isomorphic to $\C$ (see \cite [Ch. 15]{BKM}). Let $w$ be a parameter on $\hat{V}\cong\C$. Then the formulas 
\begin{equation}\label{z_{1,2}(om)}
\left\{
\begin{array}{rcl}
z_1&=&z_1(w)=-i\omega\sinh w\\\\
z_2&=&z_2(w)=-i\omega\sinh(w+i\Phi)
\end{array}
\right| w\in\C
\end{equation}
describe an infinitely sheeted covering of $\C$ onto $V$.\\
Let us ``lift'' the functions $\tilde{v}_{l}^{\beta}(z_l), z_{l}\in\C^{+}$ to $\hat{V}$. For this we must identify 
$V_l^{+}:=\Big\lbrace z\in\C^2: \Im z_l>0\Big\rbrace$ with regions on $\hat{V}$. This can be done in many ways. For example, define, for $\omega\in\C^{+}$,
\begin{equation}\label{Gamma_0'}
\Gamma_0=\Gamma_0(\omega):=\Big\lbrace w_1+i\arctan\Big(\ds\frac{\omega_1}{\omega_2}\tanh w_1\Big)\Big| w_1\in\R\Big\rbrace.
\end{equation}
Obviously, for $w\in\Gamma_0$, $\Im \big(z_1(w)\big)=0$. Moreover, 
$
\arctan\Big(\ds\frac{\omega_1}{\omega_2}\tanh w_1\Big)\xrightarrow[w_1\to\pm\infty]{}\pm\arctan\ds\frac{\omega_1}{\omega_2}.
$
For $\alpha\in\R$, define
\begin{equation*}
\Gamma_{\alpha}=\Gamma_{\alpha}(w):=\Gamma_{0}(w)+i\alpha
\end{equation*}
and for $\alpha<\beta$, define 
\begin{equation*}
V_{\alpha}^{\beta}:=\Bigg\lbrace w\in\C: \arctan\Big(\ds\frac{\omega_1}{\omega_2}\tanh w_1\Big)<\Im w<\arctan\Big(\ds\frac{\omega_1}{\omega_2}\tanh w_1\Big)+\beta\Bigg\rbrace.
\end{equation*}
For $l=1,2$, let us ``lift'' $V_{l}^{+}$ to $\hat{V}$. 
Denote this lifting  by $\hat{V}_{l}^{+}=\Big\lbrace w\in\hat{V}: \big(z_l(w)\big)\in V_{l}^{+}\Big\rbrace$. Then         
\begin{equation*}
\hat{V}_{1}^{+}=\ds\bigcup_{k=-\infty}^{\infty}{V_{2k\pi}^{(2k+1)\pi}},\quad \hat{V}_{1}^{-}=\ds\bigcup_{k=-\infty}^{\infty}{V^{2k\pi}_{(2k+1)\pi}},\quad V_{2}^{\pm}=V_{1}^{\pm}-2i\Phi.
\end{equation*}
Note that $\pm\Im \big(z_{l}(\omega, w)\big)>0,\quad w\in \hat{V}_{l}^{\pm}$.
We choose the connected component of $\hat{V}_{l}^{+}$ corresponding to the condition $\Im z_{l}>0$ as $\hat{V}_1^{+}:= V_0^{\pi}, \hat{V}_2^{+}= V_{-\Phi}^{\pi-\Phi}$,
(see Fig.\;\ref{w=i}, where $\Gamma_\alpha(w)$ are represented for $\omega_1\geq0$).

\newpage

Now we ``lift'' $\tilde{v}_{l}^{\beta}(z_{l})$ to $\hat{V}_{l}^{+},~ l=1,2,~ \beta=0,1$, using (\ref{z_{1,2}(om)}). We obtain from (\ref{tilde{v}_1^0}), (\ref{v_1^0})
\begin{align}\label{check{v}_1^0}
\hat{v}_1^0(w)=\ds\frac{i}{-i\omega\sinh w-k},\quad w\in\hat{V}_{1}^{+},\qquad \hat{v}_2^0(w)=0,\quad w\in\hat{V}_{2}^{+}.
\end{align}
Further, $\hat{v}_l^1(w)$ are analytic functions in $\hat{V}_{l}^{+}$ by (\ref{atilde{v}_l}). Our aim is to find the unknown functions $\hat{v}_{l}^1, l=1,2$. Having these functions, we obtain $\tilde{\gamma}(z)$ and the solution $v_{0}(x)$ by (\ref{8.5}), (\ref{div1}).\\ 
Note that in the case $\Phi>\pi$ the function $\tilde{\gamma}(z)$ given by (\ref{tilde{gamma}''}) is not lifted to $\hat{V}$ since $\hat{\gamma}(w):=\tilde{\gamma}\big(z_1(w), z_2(w)\big)$ is not defined at any point of $\hat{V}$. In fact, $\hat{v}_1^1(w)$ is not defined in $\hat{V}_2^{+}$, $\hat{v}_2^1(w)$ is not defined in $\hat{V}_1^{+}$ since $\hat{V}^{\ast}:=\hat{V}_{1}^{+}\cap \hat{V}_{2}^{+}=\emptyset$, see Fig.\;\ref{w=i}. In the case $\Phi<\pi$ this intersection is not empty and such lifting to $\hat{V}^{\ast}$ is possible \cite{BKM}. Thus, in that case there exists a connection between $\hat{v}_1^1$ and $\hat{v}_2^1$ generated by (\ref{div1}) since $\mathscr{H}(z)$ has zeros in $\hat{V}^{\ast}$ and $\hat{\gamma}(w)$ must vanish for $w\in \hat{V}^{\ast}$.\\
Nevertheless, a similar relation between $\hat{v}_1^1$ and $\hat{v}_2^1$  exists in the case $\Phi>\pi$ too (see \cite[{\rm Chap.\;21}]{BKM}). Let us describe  the corresponding construction. The function $\hat{\gamma}(z)$  is naturally splitted into two summands each of which is extended to $\hat{V}_{1}^{+}$ and $\hat{V}_{2}^{+}$, respectively. 
Namely
\begin{equation*}
\sin^2\Phi\;\tilde{\gamma}(z_1, z_2)=\tilde{v}_{1}(z_1, z_2)+\tilde{v}_{2}(z_1, z_2),
\end{equation*} 
where
\begin{equation*}
\left.
\begin{array}{rcl}
v_1(z_1, z_2):=\tilde{v}_1^1(z_1)+\ds\frac{z_2-2\cos\Phi\;z_1}{z_1-k},\qquad v_2(z_1,z_2):= \tilde{v}_2^1(z_2),\quad (z_1, z_2)\in\R^2.
\end{array}
\right. 
\end{equation*}
By the Paley-Wiener Theorem the function $\tilde{v}_1(z_1,z_2)$ admits an analytic continuation to $\C_{z_1}^{+}\times \C$, and $\tilde{v}_2(z_1,z_2)$ admits an analytic continuation to $\C^{+}\times\C_{z_2}^{+}$, where 
$
\C_{z_k}^{+}=\Big\lbrace z_{k}\big\vert\; \Im z_{k}>0\Big\rbrace.
$
Now we can ``lift'' $\tilde{v}_1$ and $\tilde{v}_2$ to the Riemann surface $\hat{V}$ by formulas (\ref{z_{1,2}(om)}).\\
We obtain
\begin{align}\label{check v_1}
\hat{v}_1(w)=\hat{v}_1^1(w)+\omega\sinh(w-i\Phi)\; \hat{v}_1^0(w),\quad w\in \hat{V}_1^{+},\quad\hat{v}_2(w)=\hat{v}_2^1(w), \quad w\in \hat{V}_2^{+}.
\end{align}
Then, by (\ref{check v_1}), (\ref{check{v}_1^0}),
\begin{align}\label{check{v}_1}
\hat{v}_1(w)&= \hat{v}_1^1(w)-\hat{G}(w),\quad w\in\hat{V}_1^{+},
\end{align}
where
\begin{equation}\label{G(w)}
\hat{G}(w):=\ds\frac{i \omega\sinh(w-i\Phi)}{i\omega\sinh w+k},\quad w\in\C.
\end{equation}
In the case $\Phi<\pi$, $\hat{v}_1(w)$ and $\hat{v}_2(w)$ have a common domain $\hat{V}^{\ast}$ which is not empty, and thus the connection equation has the following form:
\begin{equation}\label{cecs}
\hat{v}_1(w)+\hat{v}_2(w)=0,\quad w\in \hat{V}^{\ast}.
\end{equation} 
(see \cite[{\rm Chap.\;10}]{BKM}).\\
In the case $\Phi>\pi$ the domain $\hat{V}^{\ast}=\emptyset$ (see Fig.\;\ref{w=i}). Nevertheless, it turns out that in this case there exists a {\it connection} between $\hat{v}_1$ and $\hat{v}_2$ such that (\ref{cecs}) holds in a slightly different sense. Only in this case this equation holds for analytic continuations of $\hat{v}_1$ and $\hat{v}_2$.
Let us formulate precisely the corresponding theorem.
\begin{defn} Denote
$
\hat{V}_{\Sigma}:=\hat{V}_1^{+}\cup\hat{V}_2^{+}\cup\hat{V}^{\ast},
$
where $\hat{V}^{\ast}=V_{\pi-\Phi}^{0}$. Note that $\hat{V}_{\Sigma}=V_{-\Phi}^{\pi}$,  (see Fig.\;\ref{w=i}).
\end{defn}
\begin{teo}\label{CONEX}
(Connection equation in the case $\Phi>\pi$) \cite[Section\;20.1, Theorem\,20.1]{BKM}.
Let $v\in S'(K)$ be any distributional solution to (\ref{ra1}). Then functions (\ref{check v_1}) admit {\it analytic} continuations $[\hat{v}_{l}]$ along the Riemann surface $\hat{V}$ from $\hat{V}_{l}^{+}$ to $\hat{V}_{\Sigma}$ (see Fig.\;\ref{w=i}) and
\begin{equation}\label{CE}
\Big[\hat{v}_1(w)\Big]+\Big[\hat{v}_2(w)\Big]=0,\quad w\in\hat{V}_{\Sigma}.
\end{equation}
\end{teo}

\begin{obs}\label{r3.7}
Using the connection equation (\ref{CE}) we will find $\hat{v}_1$. The solution $u_1$ of problem (\ref{u_1}) is given by (\ref{u_1'}).
\end{obs}

\subsection{Step\;3: Reduction to a difference equation}

From (\ref{CE}), (\ref{check{v}_1}) it follows that $\hat{v}_1^1(w)$ and $\hat{v}_2^1(w)$ admit meromorphic continuations to $\hat{V}_{\Sigma}$ and 
\begin{equation}\label{al}
\hat{v}_1^1(w)+\hat{v}_2^1(w)=\hat{G}(w),\quad w\in\hat{V}_{\Sigma}.
\end{equation}
We will use the following automorphisms on $\hat{V}$ $\Big({\rm see}~ \cite[{\rm Ch.}\;13]{BKM} ~{\rm and}~ \cite[(73)]{KMM}\Big)$:
\begin{equation}\label{h_1h_2}
\left.
\begin{array}{rcl}
h_1 w= -w+\pi i, \qquad h_2 w=-w+\pi i-2i\Phi, \quad w\in\C
\end{array}
\right. 
\end{equation}
which are symmetries with respect to $i\ds\frac{\pi}{2}$ and $i\ds\frac{\pi}{2}-i\Phi$, respectively.\\
Sometimes we will use the notation $f^{h_{l}}(w):=f\Big(h_{l}(w)\Big),~l=1,2$.\\
The functions $\hat{v}_1^1$ and $\hat{v}_2^1$ are automorphic functions with respect to $h_1$ and $h_2$, respectively: 
\begin{align}\label{h_1}
&\hat{v}_1^{1h_{1}}(w)=\hat{v}_1^1(-w+\pi i)=\hat{v}_1^1(w),\quad w\in\hat{V}_1^{+},
\\ \nonumber\\\label{h_2}
&\hat{v}_2^{1h_{2}}(w)=\hat{v}_2^1(-w+\pi i-2i\Phi)=\hat{v}_2^1(w),\quad w\in\hat{V}_2^{+},
\end{align}
as follows from the fact that $\tilde{v}_{l}^1(z_l)$ depend only on $z_l$ and hence their liftings $\hat{v}_{l}(w)$ to $\hat{V}_{l}^{+}$ satisfies (\ref{h_1}), (\ref{h_2}) since $\sinh w$ satisfy (\ref{h_1}) and $\sinh(w+i\Phi)$ satisfies (\ref{h_2}).\\
Thanks to this automorphy we can eliminate one unknown function in the undetermined equation (\ref{al}) and reduce it to an equation with a shift, see \cite{KMM}. The idea of this method is due to Malyshev \cite{Mal}.
\begin{lem}\label{ll-1}
Let $v\in S'(\R)$ satisfy (\ref{ra1}-\ref{ra3}) and $v_l^1(x_l), l=1,2$, be its Neumann data. Then their liftings to $\hat{V}$, $\hat{v}_{l}^1(w), w\in\hat{V}_{l}^{+}$, admit meromorphic continuations to $\C$ (which we also denote $\hat{v}_{l}^1$) such that for $w\in\C$
\begin{align}\label{3.25}
\hat{v}_1^1(w)+\hat{v}_2^1(w)&=\hat{G}(w),
\end{align}
and they are $h_l$-automorphic functions,
\begin{align}\label{icheck{v}_1^1}
\hat{v}_1^1\big(h_1 (w)\big) &= \hat{v}_1^1(w),
\end{align}
$$
\hat{v}_2^1\big(h_2 (w)\big) = \hat{v}_2^1(-w+\pi i-2i\Phi)= \hat{v}_2^1(w).
 $$
\end{lem}
The proof of this lemma is given in Appendix\;\ref{aal}.

\medskip
\no Now we reduce system (\ref{3.25})-(\ref{icheck{v}_1^1}) to a difference equation, which is also called the shift equation. This reduction is the part of MAF which was introduced in \cite{Mal} for difference equations in angles. It uses the automorphy of $\hat{v}_{l}^{\beta}$ on $\hat{V}$ under the automorphisms $h_l$ and the term MAF is due to this observation.\\
Define, for $w\in\C$, 
\begin{equation}\label{d21}
\hat{G}_2(w):=\hat{G}(w)-\hat{G}\Big(h_2(w)\Big)=\ds\frac{i\omega\sinh(w-i\Phi)}{i\omega\sinh w+k}-\ds\frac{i\omega\sinh(w+3i\Phi)}{i\omega\sinh(w+2i\Phi)+k}.
\end{equation}
For a region $U$ in $\C$ we will denote  here and everywhere below $\mathcal{M}(U)$ the set of meromorphic functions on $U$.
\begin{lem}\label{Vin}
Let $v\in S'(K)$ satisfy (\ref{ra1})-(\ref{ra2}). Then the connection equation (\ref{CE}) holds, the function $\hat{v}_1^1$ belongs to $\mathcal{M}(\C)\cap \mathcal{H}(\hat{V}_1^{+})$,
satisfies the difference equation  
\begin{equation}\label{de.}
\begin{array}{lll}
\hat{v}_1^1(w)-\hat{v}_1^1(w+2i\Phi)=\hat{G}_2(w),
\end{array}
\end{equation}
 and the automorphic condition (\ref{icheck{v}_1^1}).
\end{lem}
The proof of this lemma is given in Appendix\;\ref{ade}.

\medskip

Our goal is to find $\hat{v}_1^1(w)\in \mathcal{M}(\C)\cap \mathcal{H}\big(\hat{V}_1^{+}\big)$ s.t. (\ref{de.}), (\ref{icheck{v}_1^1}) and the condition 
$
\hat{v}_1(w)\in \mathcal{H}\big(\hat{V}_{\Sigma}\big) 
$ 
hold. Here $\hat{v}_1$ is given by (\ref{check v_1}). In its turn, this condition is equivalent to the condition
\begin{equation}\label{Si2}
\hat{v}_1(w)=\hat{v}_1^1(w)-\hat{G}(w)\in\mathcal{H}\big(\hat{V}_{\Sigma}\big)
\end{equation} 
by (\ref{check{v}_1}).\\ 
In the next section we find the necessary and sufficient conditions for $\hat{v}_1^1$ such that condition (\ref{Si2}) holds.

%\newpage
\subsection{Necessary and sufficient condition for $\hat{v}_1^1$}

The analyticity condition (\ref{Si2}), which follows from the connection equation (\ref{CE}), imposes certain necessary conditions for the poles of the function $\hat{v}_1^1$, more exactly,  of  its continuation obtained in Lemma\;\ref{ll-1}. This section is devoted to the derivation of these conditions and to the proof of the fact that they are also sufficient for (\ref{CE}) to hold.\\
We will often use the following evident statement.   
\begin{lem}\label{lh_1}
Let $A\in\mathcal{M}(\C)$, and $A$ satisfy (\ref{icheck{v}_1^1}). Then
$
\overset{}{\underset{w=w_0}{res}}\;A(w)=-\overset{}{\underset{w=-w_0+\pi i}{res}}\;A.
$
\end{lem}

\newpage
\begin{figure}[htbp]
\centering
\includegraphics[scale=0.33]{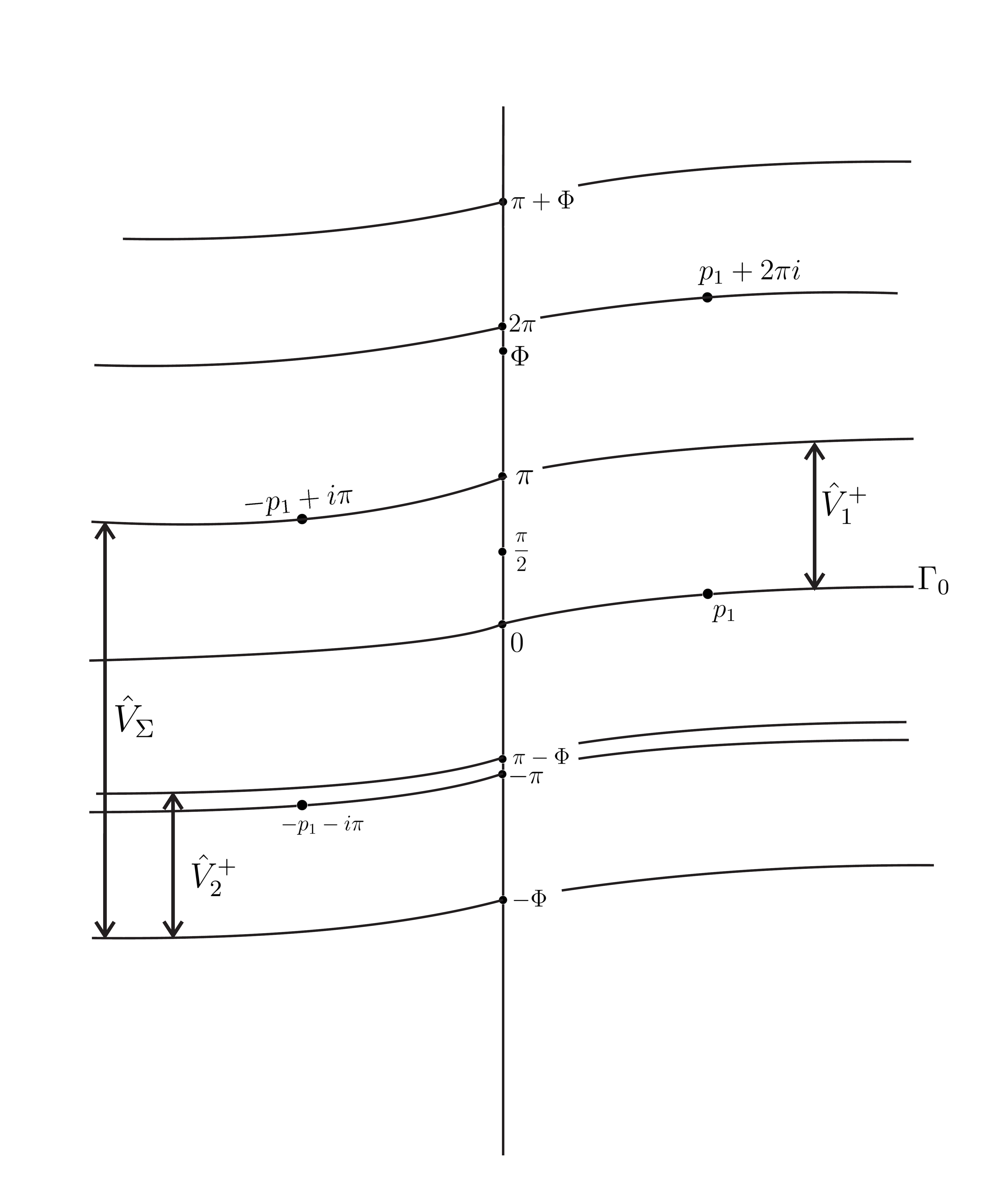}
\caption{Necessary conditions}\label{AJ}
\end{figure}

Denote
\begin{equation}\label{P}
P:=\Big\lbrace p_1, -p_1\pm \pi i, p_1+2\pi i\Big\rbrace,\qquad p_1:=\sinh^{-1}\Big(\ds\frac{ik}{\omega}\Big)\in\Gamma_0.
\end{equation}
(See Fig.\;\ref{AJ}, where the positions  of the curves  $\Gamma_\alpha$ correspond to the case $\Re \omega>0$.  We will always assume in the following that this is the case; in the case $\Re \omega <0$ the construction is similar.)
Introduce the following two important parameters:
\begin{equation}\label{r1r2}
r_1:=\overset{}{\underset{-p_1\pm \pi i}{res}}\;\hat{G}:=-\ds\frac{\sinh(p_1+i\Phi)}{\cosh p_1},\quad r_2:=\overset{}{\underset{p_1}{res}}\;\hat{G}:=\ds\frac{\sinh(p_1-i\Phi)}{\cosh p_1}.
\end{equation}
In the next  proposition we give a necessary and sufficient condition for $\hat{v}_1^1$ guaranteeing that  condition (\ref{d21}) holds.
\begin{prop}\label{pns}
Let $\hat{v}_1^1\in\mathcal{M}(\C)\cap \mathcal{H}(\hat{V}_{1}^{+})$ satisfy (\ref{de.}) and (\ref{icheck{v}_1^1}). Then $\hat{v}_1\in\mathcal{H}\big(\hat{V}_{\Sigma}\big)$ if and only if 
\begin{equation}\label{cu1}
\hat{v}_1^1\in\mathcal{H}\big(V_{-\Phi}^{\Phi+\pi}\setminus P\big)
\end{equation}
and
\begin{equation}\label{u2}
\overset{}{\underset{p_1}{res}}\;\hat{v}_1^1=r_2,\quad \overset{}{\underset{-p_1-\pi i}{res}}\;\hat{v}_1^1=r_1.
\end{equation}
\end{prop}
{\bf Proof}. First, let us prove the necessity of condition (\ref{u2}). By (\ref{G(w)}) the poles of $\hat{G}(w)$ in $\C$ are
$
p_1+2k\pi i, -p_1-\pi i+2k\pi i, k\in\Z.
$
Of all of these poles only $p_1, -p_1-\pi i$ belong to $\hat{V}_{\Sigma}$ (see Fig.\;\ref{AJ}). 
Hence formulas (\ref{u2}) follow from (\ref{Si2}), (\ref{r1r2}). Further, since $\hat{v}_1^1\in \mathcal{H}\Big(V_{-\Phi}^{\pi}\setminus \lbrace p_1, -p_1-\pi i\rbrace\Big)$ by (\ref{Si2}), $\hat{v}_1^1\in \mathcal{H}\Big(V_{0}^{\pi+\Phi}\setminus \lbrace p_1+2\pi i, -p_1+\pi i\rbrace\Big)$ by (\ref{icheck{v}_1^1}). Hence (\ref{cu1}) also holds.\\
Let us prove the sufficiency of conditions (\ref{cu1}), (\ref{u2}). From (\ref{u2}) and (\ref{r1r2}) it follows that
$
\overset{}{\underset{w=p_1}{res}}\;\hat{v}_1(w)=\overset{}{\underset{w=-p_1-\pi i}{res}}\;\hat{v}_1(w)=0.
$ 
Hence, $\hat{v}_1\in\mathcal{H}\big(\hat{V}_{\Sigma}\big)$ since $\hat{v}_1^1\in\mathcal{H}\Big(\hat{V}_{\Sigma}\setminus\lbrace p_1, -p_1-\pi i\rbrace\Big)$ by (\ref{cu1}) and $\hat{G}(w)$ belongs to this space too. The proposition is proven.~~~~~$\blacksquare$

\begin{obs}
Condition (\ref{cu1}) implies $\hat{v}_{1}^1\in\mathcal{H}\big(\hat{V}_1^{+}\big)$.
\end{obs}

\section{$h_1$-automorphic solution of difference equation (\ref{de.}), $\Phi\neq\ds\frac{3\pi}{2}$}
\setcounter{equation}{0}

In this section we construct an $h_1$-automorphic solution of difference equation (\ref{de.}) satisfying all conditions of Proposition\;\ref{pns} for $\Phi\neq\ds\frac{3\pi}{2}$. 
This limitation is related to the method of obtaining a solution which uses the Cauchy-type integral. The kernel of this integral must be analytic on the integration contour. In turns out that it is possible to find such a kernel only when $\Phi\neq\ds\frac{3\pi}{2}$. Fortunately, the case $\Phi=\ds\frac{3\pi}{2}$ does not need an integral of the Cauchy-type since the difference equation (\ref{de.}) is solved by elementary methods in this case (see Section\;5).

\subsection{Poles of $\hat{G}_2$ and asymptotics.}
In this subsection we give the properties of $G_2$ which are necessary for the solution of the main problem.
Let 
\begin{align}\label{pq}
 \mathscr{P}:=\Big\lbrace p_1+2\pi i k-2i\Phi m: k\in\Z, m=0,1\Big\rbrace,\quad
\mathscr{Q}:=\Big\lbrace -p_1+\pi i+2\pi i k-2i\Phi m: k\in\Z, m=0,1\Big\rbrace.
\end{align}
Obviously, the poles of $\hat{G}_2$ belong to $\mathscr{P}\cup\mathscr{Q}$ by (\ref{d21}).\\
Denote 
$
q_1=-p_1-\pi i+2i\Phi.
$
\begin{lem}\label{lpG2}
i) The poles of $\hat{G}_2$ belonging to $\overline{V_{\frac{\pi}{2}-\Phi}^{\pi-\Phi}}$ are 
$-p_1-\pi i~{\rm for}~\Phi\geq\frac{3\pi}{2} ~{\rm and}~ -q_1+\pi i~{\rm for}~\Phi\leq\frac{3\pi}{2}$.\\
The residues of $\hat{G}_2$ at these points are 
\begin{equation*}
\overset{}{\underset{-p_1-\pi i}{res}}\;\hat{G}_2= \overset{}{\underset{-q_1+\pi i}{res}}\;\hat{G}_2=r_1.
\end{equation*}
ii) The poles of $\hat{G}_2$ belonging to $\overline{V_{-\frac{\pi}{2}-\Phi}^{-\Phi}}$ are $p_1-2\pi i$ for $\Phi\geq\ds\frac{3\pi}{2}$, $-p_1+\pi i-2i\Phi$ for $\Phi\leq \ds\frac{3\pi}{2}$ and
\begin{equation*}
\overset{}{\underset{p_1-2\pi i}{res}}\;\hat{G}_2=\overset{}{\underset{-p_1+\pi i-2i\Phi}{res}}\;\hat{G}_2=r_2.
\end{equation*}  
iii) The function $\hat{G}_2$ admits the asymptotics 
\begin{equation}\label{asG2}
\hat{G}_2(w)=\mp 2i\sin\Phi\Big(1+O(e^{\mp \Re w})\Big),\quad \Re w\to \pm\infty
\end{equation}
uniformly with respect to $\Im w$.\\ 
iv)
\begin{equation}\label{Shat{G}}
\hat{G}_2\Big(h_2(w)\Big)=-\hat{G}_2(w),\quad\hat{G}_2\Big(\ds\frac{\pi}{2}-i\Phi\Big)=0.
\end{equation}
\end{lem}
{\bf Proof}. The first three  assertions follow  directly from (\ref{d21}). The last assertion (\ref{Shat{G}}) follows from the fact that $h_2\Big(\ds\frac{\pi i}{2}-i\Phi\Big)=\ds\frac{\pi i}{2}-i\Phi$. ~~~~$\blacksquare$
\begin{obs}\label{ras}
In the case $\Phi=\ds\frac{3\pi}{2}$, 
$$
\hat{G}_2(w)=\pm 2i\Big(1+O(e^{\mp 2\Re w})\Big),\quad \Re w\to\pm\infty.
$$
However, this will not affect the final results.
\end{obs}

\subsection{Reduction of problem (\ref{de.}), (\ref{icheck{v}_1^1}) to a conjugate problem}

\begin{figure}[htbp]
\centering
\includegraphics[scale=0.33]{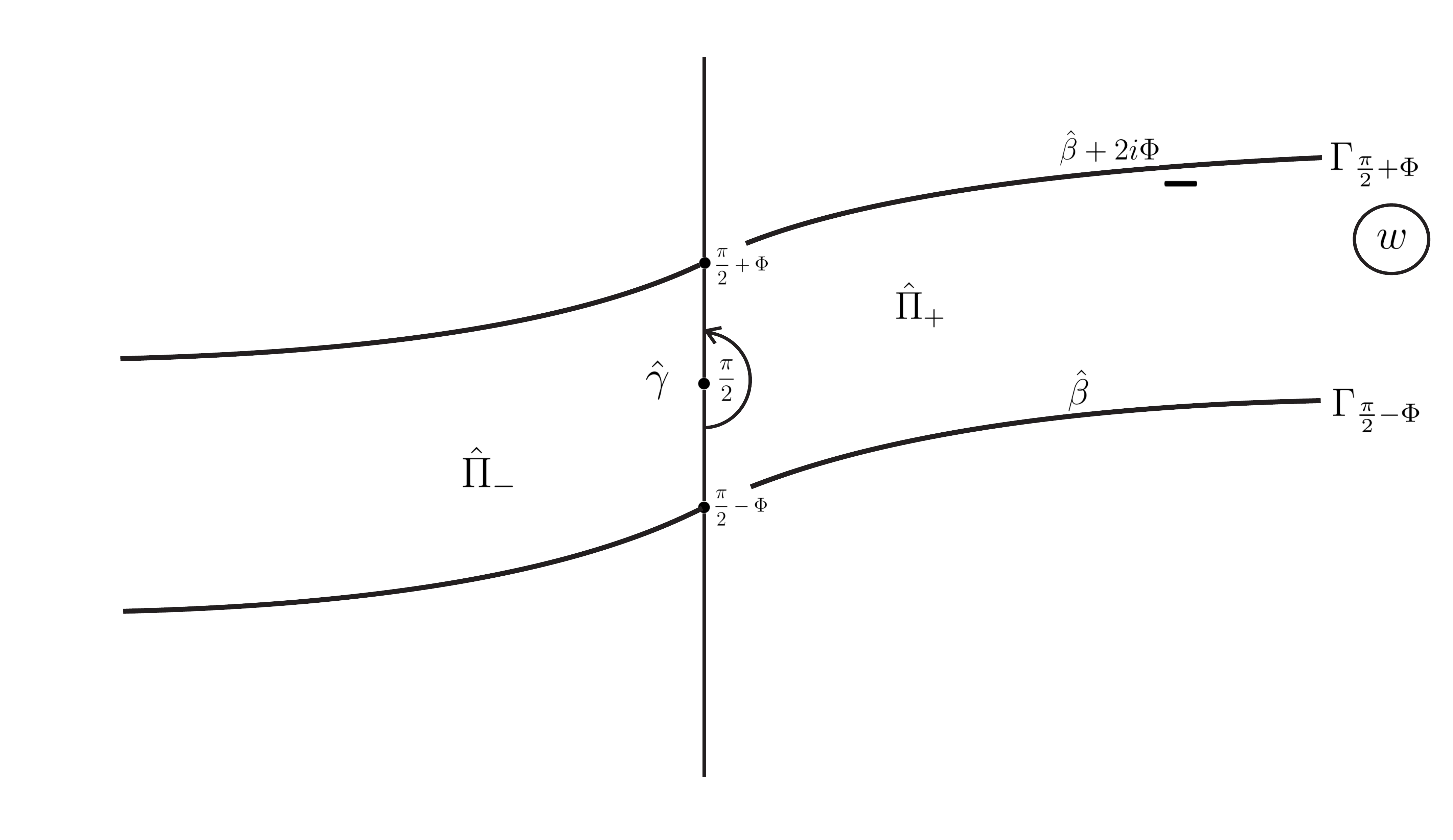}
\caption{Reduction to conjugate problem}\label{APi}
\end{figure}

Denote
$$
\hat{\Pi}:=V_{\frac{\pi}{2}-\Phi}^{\frac{\pi}{2}+\Phi},\qquad \hat{\Pi}_{\pm}=\Big\lbrace w\in\hat{\Pi}:\pm\Re w>0\Big\rbrace,\quad \partial\hat{\Pi}_{+}=\hat{\beta}\cup\hat{\gamma}\cup(\hat{\beta}-2i\Phi),
$$ 
where
\begin{equation*}
\hat{\beta}:=\Big\lbrace w\in\Gamma_{\frac{\pi}{2}-i\Phi}: \Re w\geq0\Big\rbrace,\quad \hat{\gamma}:=\Big\lbrace w\big\vert \Re w=0,\quad \Im w\in\Big[\frac{\pi}{2}-\Phi, \frac{\pi}{2}+\Phi\Big]\Big\rbrace.
\end{equation*}
We will look for a solution of the following problem: to find an analytic function in $\hat{\Pi}_{+}$ whose  boundary values on $\partial\hat{\Pi}_{+}$,
\begin{align*}
&\hat{a}_{1}(w+i0),\quad w\in \hat{\beta};\quad\hat{a}_{1}(w+2i\Phi-i0),\quad w\in \hat{\beta}+2i\Phi, \quad\hat{a}_{1}(w+0),\quad w\in\hat{\gamma},
\end{align*}
exist and are such that they satisfy the following conditions of conjugation
\begin{align}\label{C1}
&\hat{a}_{1}(w+i0)-\hat{a}_{1}(w+2i\Phi-i0)=\hat{G}_2(w),\quad w\in\hat{\beta}
\\\nonumber\\\label{C2}
&\hat{a}_{1}(w+0)=\hat{a}_{1}(-w+\pi i+0),\quad w\in\hat{\gamma},
\end{align}
see Fig.\;\ref{APi}

\newpage

\subsection{Solution of the conjugate problem (\ref{C1}), (\ref{C2}) for $\Phi\neq\ds\frac{3\pi}{2}$}
\begin{figure}[htbp]
\centering
\includegraphics[scale=0.33]{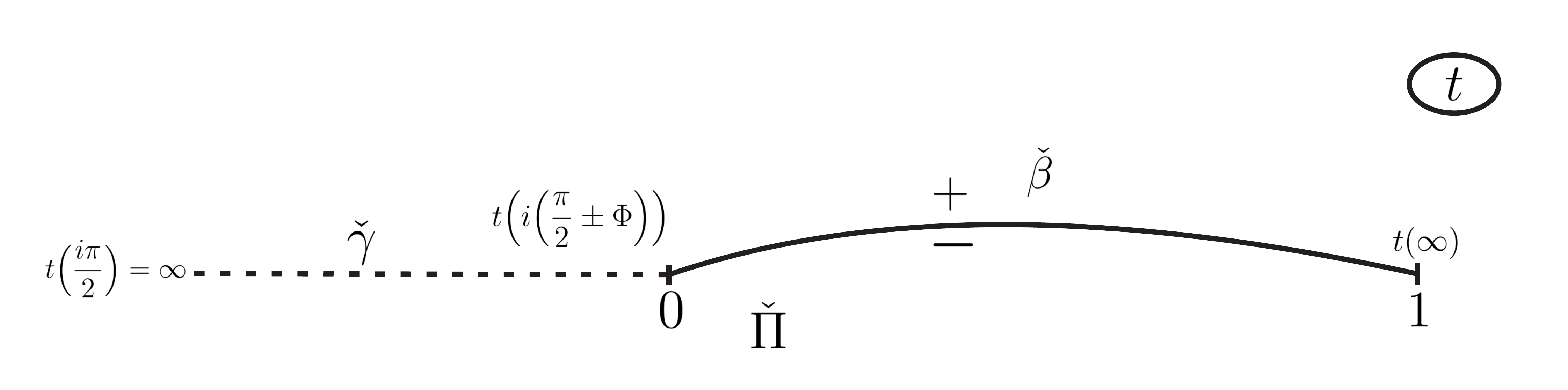}
\caption{Riemann-Hilbert problem}\label{At}
\end{figure}

We start solving the problem (\ref{C1}), (\ref{C2}). We will reduce this problem to a Riemann-Hilbert problem. To this end we map $\hat{\Pi}_{+}$ conformally onto $\check{\Pi}:=\C^{\ast}\setminus \check{\beta}$, where $\C^{\ast}$ is the Riemann sphere, $\check{\beta}:=t(\beta)$. For example, define
\begin{equation}\label{t(w)}
w\mapsto t=t(w)=\coth^2\Big(\ds\frac{\pi}{2\Phi}\Big(w-\ds\frac{\pi i}{2}\Big)\Big),\quad w\in\hat{\Pi}_{+}. 
\end{equation}  
Denote the inverse transform $w(t): \check{\Pi}_{+}$ to $\hat{\Pi}_{+}$.\\
Note that when $w\in\hat{\Pi}_{+}$  tends to $\hat{\beta}$ implies that 
 $t(w)$ tends  to $\check{\beta}$ from above, and when $w\in\hat{\Pi}_{+}$
tends to $\beta+2i\Phi$ implies that  $t(w)$  to $\check{\beta}$ from below.   
Obviously,
\begin{align*}
t(\hat{\gamma})=:\check{\gamma}=[-\infty,0],\quad t\Big(\pm i\Big(\ds\frac{\pi}{2}+\Phi\Big)\Big)=0,\quad t(\infty)=1,\quad t\Big(\ds\frac{\pi i}{2}\Big)=\infty
\end{align*}
 (see Fig.\ref{At}).

\vspace{0.5cm}
Then problem (\ref{C1}), (\ref{C2}) is equivalent to the Riemann-Hilbert problem for $\check{a}_1(t):=\hat{a}_{1}\big(w(t)\big), t\in\check{\Pi}_{+}$, which at the same time is the saltus problem (see \cite[Ch. 16, 18]{BKM})
\begin{equation}\label{RG} 
\check{a}_{1}(t+i0)-\check{a}_{1}(t-i0)=\check{G}_2(t),\quad t\in\check{\beta}.
\end{equation}
Here $\check{G}_2(t)=\hat{G}_2\Big(w(t)\Big), t\in\hat{\Pi}$, $\check{a}_{1}(t):=\check{a}_{1}\big(w(t)\big),  t\in\check{\Pi}_{+}$;  $\check{a}_{1}(t\pm i0)=\ds\lim_{\varepsilon\to 0+} \check{a}_{1}(t\pm i\varepsilon)$ and for $t\in\check{\beta}$,  $w\in\hat{\beta}$, $w\in\overline{\hat{\Pi}_{+}}$. From (\ref{t(w)}) and (\ref{Shat{G}}) it follows that $\check{G}_2(t)$ and $\check{G}'_2(t)$ are continuous on $\overline{\check{\beta}}$ and 
\begin{equation}\label{check{G}_2(0)=0}
\check{G}_2(0)=0,\qquad \check{G}_2(1)=-2i\sin\Phi.
\end{equation}
It is well known that a particular  solution of (\ref{RG}) is given by the Cauchy type integral
\begin{align}\label{cons}
\check{a}_{1}(t)=\ds\frac{1}{2\pi i}\ds\int\limits_{\check{\beta}}\ds\frac{\check{G}_2(t')}{t'-t}\;dt',\quad t\in\check{\Pi}.
\end{align}
Obviously $\check{a}_{1}(t)\in\mathcal{H}(\check{\Pi})$, and $\check{a}_{1}(t)\xrightarrow[t\to \infty]{}0$. Moreover, by (\ref{check{G}_2(0)=0}) there exists $\ds\lim_{t\to 0}\;\check{a}_1(t), t\notin\check{\beta}$.\\
In the following lemma we establish an  asymptotics of  (\ref{cons})
at $t=1$;  it plays an important role in describing the fact that the 
solution belongs to a certain class and hence  its uniqueness.
\begin{lem}\label{l as check{a}_1}
The function $\check{a}_1(t)$ admits the following asymptotics
\begin{align}\label{as hat{a}_1}
\check{a}_{1}(t)=-\ds\frac{\check{G}_2(1)}{2\pi i}\ln\ds\frac{1}{t-1}+C+O(t-1),\quad t\to 1, 
\end{align}
where, by  $\ln\ds\frac{1}{t-1}$, we understand a certain branch that is single-valued on a plane cut along $\check{\beta}$ and $C$ depends only on $\hat{G}_2$; moreover,
\begin{equation}\label{4.22'}
\ds\frac{d}{dt}\check{a}_1(t)=\ds\frac{\sin\Phi}{\pi}\;\ds\frac{1}{1-t}+\ds\frac{\check{G}'_2(1)}{2\pi i} \ln(1-t)+C_1+o(1),
\end{equation}
where $C$, $C_1$ depend only on $\check{G}_2$.
\end{lem}
{\bf Proof.} 
(\ref{as hat{a}_1}) follows from (\ref{cons}) (see \cite[\S 16]{Mus}). Let us find the asymptotics of $\ds\frac{d}{dt}\check{a}_1(t),~t\to1$. (\ref{cons}) gives
\begin{align*}
\ds\frac{d}{dt} \check{a}_1(t)=\ds\frac{1}{2\pi i}\ds\int\limits_{\check{\beta}}\ds\frac{\check{G}_2(t') dt'}{(t'-t)^2}.
\end{align*}
Represent $\check{G}_2(t')$ in the form $\check{G}_2(t')=\check{G}_{2}(1)+\check{G}_{2}(1)(t'-1)+\psi(t')(t'-1)^2, ~ t'\in\check{\beta}$, where $\psi(t')\in C^{\infty}\big(\overline{\check{\beta}}\big)$. This is possible since $\check{G}_{2}(t')\in C^{\infty}\big(\overline{\check{\beta}}\big)$ by (\ref{pq}) if $\Phi\neq\ds\frac{3\pi}{2}$.
Then 
\begin{equation*}
\begin{array}{lll}
\ds\frac{d}{dt}\check{a}_1(t)
=\ds\int\limits_{\check{\beta}}\ds\frac{\psi(t')(t'-1)^2}{(t'-t)^2} dt'+\ds\frac{1}{2\pi i}\ds\int\limits_{\check{\beta}}\ds\frac{\check{G}_2(1)+\check{G}_2'(1)(t'-1)}{(t'-t)^2} dt'.
\end{array}
\end{equation*}
Obviously 
$
\ds\frac{1}{2\pi i}\ds\int\limits_{\check{\beta}}\ds\frac{(t'-1)^2\;\psi(t')}{(t'-1)^2} dt'=C_1+o(1),\quad t\to1,
$
where $C_1$ depends only on $\check{G}_2$;
$
\ds\frac{1}{2\pi i}\ds\int\limits_{\check{\beta}}\ds\frac{\check{G}_2(1)}{(t'-t)^2} dt'
= \ds\frac{\sin\Phi}{\pi} \ds\frac{1}{1-t}-\ds\frac{\sin\Phi}{\pi}+o(1),\quad t\to1
$
by (\ref{check{G}_2(0)=0}). Further,
$
\ds\frac{1}{2\pi i}\ds\int\limits_{\check{\beta}}\ds\frac{\check{G}'_2(1) (t'-1)}{(t'-t)^2} dt'
= \ds\frac{\check{G}'_2(1)}{2\pi i}\Big[\ln(1-t)\Big]+C_2+o(1),\quad t\to 1.
$
This implies (\ref{4.22'}).~~~~~$\blacksquare$

\begin{figure}[htbp]
\centering
\includegraphics[scale=0.27]{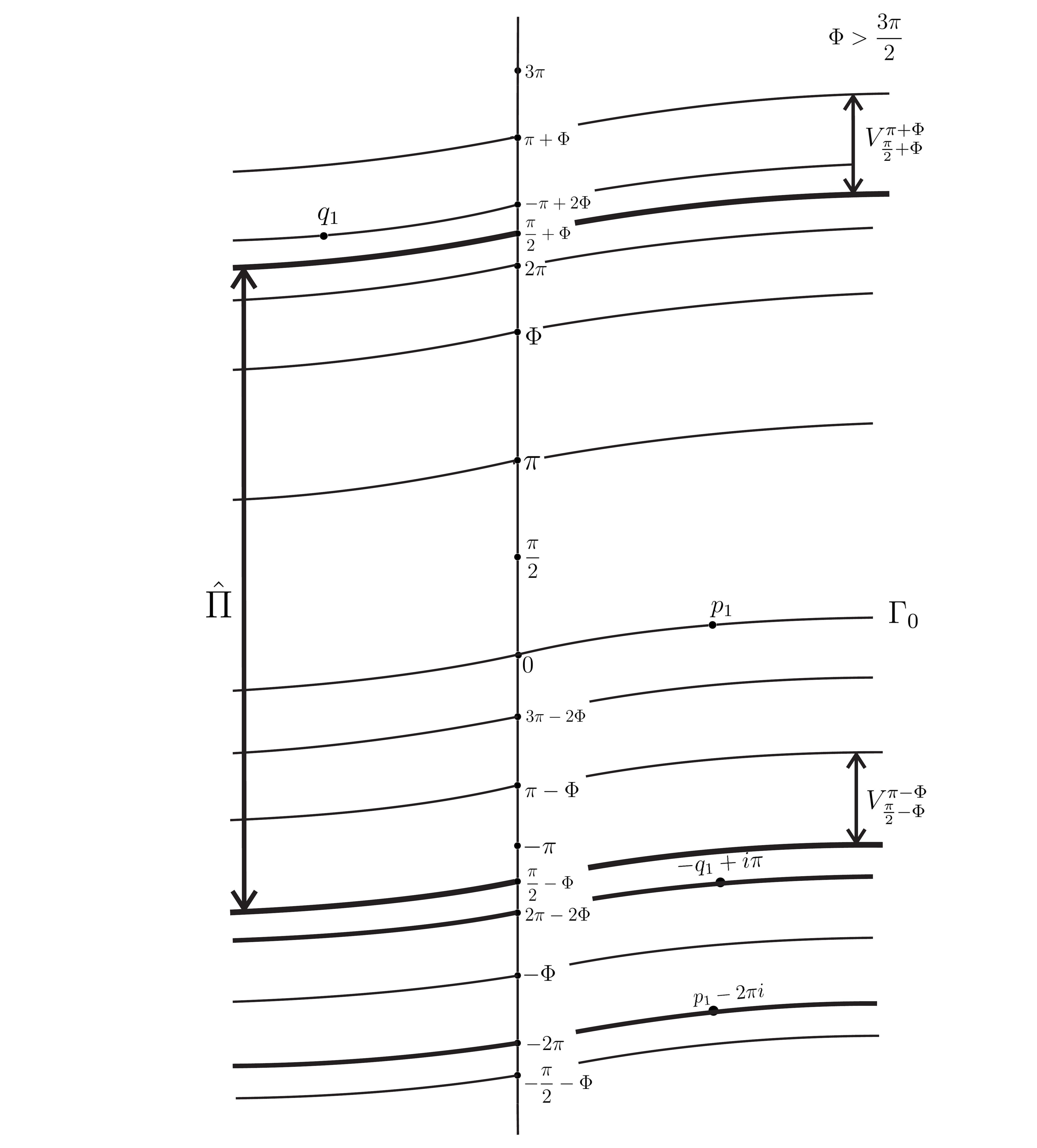}
\caption{$\Phi>\ds\frac{3\pi}{2}$}\label{phi+}
\end{figure}

\begin{figure}[htbp]
\centering
\includegraphics[scale=0.33]{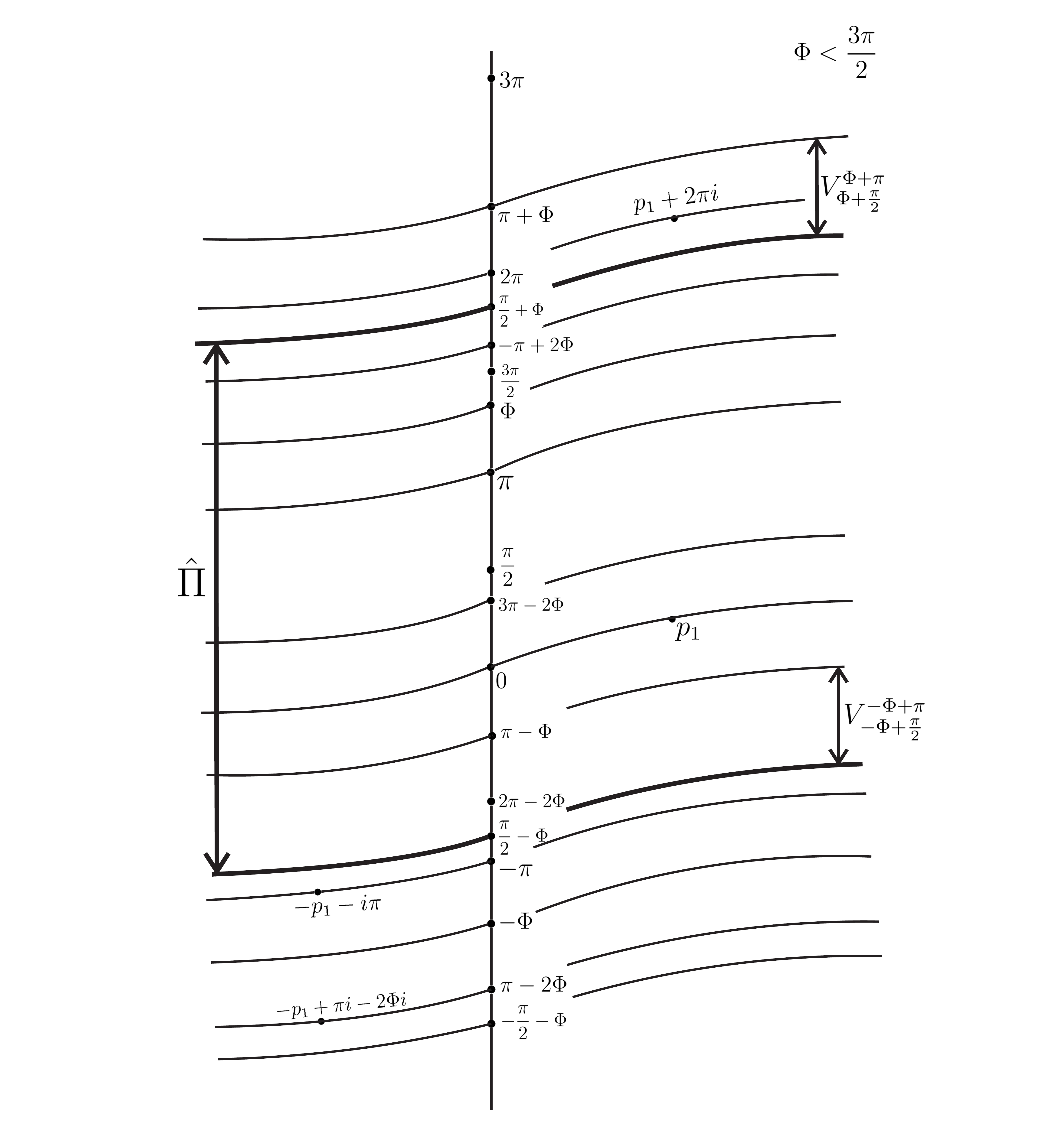}
\caption{$\Phi<\ds\frac{3\pi}{2}$}\label{Phi+1}
\end{figure}

\medskip
Now we are able to find a solution to problem (\ref{de.}). First we define this solution in $\overline{\hat{\Pi}}$ and then we extend it to $\C$.\\
Let us define $\hat{a}_1(w), w\in\overline{V_{-\frac{\pi}{2}-\Phi}^{-\Phi}}$ by the formula 
\begin{align}\label{hat{a}_1(w)}
& \hat{a}_{1}(w):=\check{a}_1\big(t(w)\big)+C_2, \quad w\in\hat{\Pi}_{+},
\end{align}
where $\check{a}_1(t)$ is given by (\ref{cons}) and
\begin{equation}\label{defC_1}
C_2=\ln 4\;\ds\frac{\sin\Phi}{\Phi}-C,
\end{equation}
with $C$ taken from (\ref{as hat{a}_1}).\\
Obviously $\hat{a}_1(w)$ satisfies (\ref{C1}) and (\ref{C2}), since a constant is a solution of the homogeneous equation (\ref{C1}) and satisfies (\ref{C2}). 
Moreover, the same formula (\ref{hat{a}_1(w)}) defines the analytic function $\hat{a}_1(w)$ in $\overline{\hat{\Pi}}$, since $\hat{a}_1(w)$ satisfies (\ref{C2}). Obviously, (\ref{t(w)}) implies that
\begin{align}\label{au}
& \hat{a}_{1}(w)=\hat{a}_1(-w+\pi i), \quad w\in\overline{\hat{\Pi}}.
\end{align}
Moreover, $\hat{a}_1(w)$ satisfies (\ref{C1}).
In fact, for $w\in\hat{\beta}$, $\Re w>0$ this follows from (\ref{RG}), and for $w\in\hat{\beta}$, $\Re w<0$ this follows from (\ref{RG}) and (\ref{Shat{G}}).
Further, we extend  $\hat{a}_{1}(w)$ to $\overline{{V}_{\frac{\pi}{2}-3\Phi}^{\frac{\pi}{2}+3\Phi}}$ by formulas corresponding to the difference equation (\ref{de.}):
\begin{equation}\label{e1}
\hat{a}_{1}(w)=\hat{a}_{1}(w-2i\Phi)-\hat{G}_2(w-2i\Phi),\quad w\in \overline{\hat{\Pi}}+2i\Phi,\quad \overline{\hat{\Pi}}+4i\Phi,\cdots
\end{equation}
and 
\begin{equation}\label{3.4'}
\hat{a}_{1}(w)=\hat{a}_{1}(w+2i\Phi)+\hat{G}_2(w),~~{\rm for}~~ w\in\overline{\hat{\Pi}}-2i\Phi,\quad \overline{\hat{\Pi}}-4i\Phi,\cdots
\end{equation}
This extension is meromorphic by (\ref{C1}) and still has property (\ref{au}). Let us  prove this.
Let $w\in\hat{\Pi}-2i\Phi$ (see Fig.\;\ref{phi+}, \ref{Phi+1}). Then $w+\pi i\in\hat{\Pi}+2i\Phi$. By (\ref{3.4'}) and (\ref{e1}) we have $\hat{a}_1(w)=\hat{a}_1(w+2i\Phi)+\hat{G}_2(w), ~ \hat{a}_1(-w+\pi i)=\hat{a}_1(-w+\pi i-2i\Phi)-\hat{G}_2(-w+\pi i-2i\Phi)$. But $\hat{a}_1(w+2i\Phi)=\hat{a}_1(-w+\pi i-2i\Phi)$ since $\hat{a}_1\Big(h_1(w)\Big)=\hat{a}_1(w),~w\in\overline{\hat{\Pi}}$ and $\hat{G}_2(w)=-\hat{G}_2(-w+\pi i-2i\Phi)$ by (\ref{Shat{G}}). Hence $\hat{a}_1(w)=\hat{a}_1(-w+\pi i)$ in this case. Similarly this is true for $w\in\hat{\Pi}+2i\Phi$.~~~~~~$\blacksquare$
\begin{obs}
Similarly, $\hat{a}_1(w)$ admits a meromorphic extension to $\C$ which satisfies (\ref{de.}) and (\ref{icheck{v}_1^1}).
\end{obs}
\begin{prop}\label{P_2a_2}
For $\Phi\neq\ds\frac{3\pi}{2}$\\
i) There exists a meromorphic in $\C$ and analytic in $\hat{\Pi}$ solution $\hat{a}_{1}$ of problem (\ref{de.}), (\ref{icheck{v}_1^1}) given by (\ref{e1}), (\ref{3.4'}).\\ 
ii) The function $\hat{a}_1(w)$ admits the following asymptotics 
\begin{align}\label{ashat{a}_1}
& \hat{a}_1(w)=\pm\ds\frac{\sin\Phi}{\Phi}\big(w-\ds\frac{\pi i}{2}\Big)+o\Big(e^{\mp\frac{\pi}{2\Phi}w}\Big),\quad \Re w\to\pm\infty,
\\\nonumber\\\label{4.27'}
& \ds\frac{d}{dw}\hat{a}_1(w)=\pm\ds\frac{\sin\Phi}{\Phi}+o\big(e^{\mp \frac{\pi}{\Phi}w}\big),\quad \Re w\to\pm\infty,
\end{align}
uniformly with respect to $\Im w, w\in\hat{\Pi}_{+}$.\\
iii) The solution $\hat{a}_1$ has poles in $\overline{V_{-\frac{\pi}{2}-\Phi}^{\pi+\Phi}}$ for $\Phi>\ds\frac{3\pi}{2}$ only at $q_1:=-p_1-\pi i+2i\Phi,~ -q_1+\pi i=p_1+2\pi i-2i\Phi$ and  $p_1-2\pi i$ with residues
\begin{equation}\label{rcheck_1^1}
\overset{}{\underset{q_1}{res}}\;\hat{a}_{1}=-r_1,\qquad \overset{}{\underset{-q_1+\pi i}{res}}\;\hat{a}_{1}=r_1,\quad \overset{}{\underset{p_1-2\pi i}{res}}\;\hat{a}_{1}=r_2.
\end{equation}
For $\Phi\leq\ds\frac{3\pi}{2}$, $\hat{a}_1$ has poles in $\overline{V_{-\frac{\pi}{2}-\Phi}^{\pi+\Phi}}$ only at $p_1+2\pi i$,  $-p_1-\pi i$ and $-p_1+\pi i-2i\Phi$ with residues 
\begin{equation}\label{3.4''}
\overset{}{\underset{p_1+2\pi i}{res}}\;\hat{a}_1=-r_1,\qquad \overset{}{\underset{-p_1-\pi i}{res}}\;\hat{a}_1=r_1,\quad \overset{}{\underset{-p_1+\pi i-2i\Phi}{res}}\;\hat{a}_1=r_2.
\end{equation}
\end{prop}

{\bf Proof.} Statement $i)$ is proved above. The asymptotics (\ref{ashat{a}_1}), (\ref{4.27'}) are proved in Appendix\;\ref{AC1}.\\
Statement $iii)$ follow from the difference equation (\ref{de.}), $h_1$-automorphicity of $\hat{a}_1$, Lemma\;\ref{lh_1},  Lemma\;\ref{lpG2},
since the function $\hat{a}_1$ is analytic in $\hat{\Pi}$ by i), see Fig.\;\ref{phi+}, \ref{Phi+1}.~~~~~$\blacksquare$

\medskip
As we will see below this asymptotics coincides with the asymptotics of the function $\hat{v}_1^1$ in the case $\Phi=\ds\frac{3\pi}{2}$.
In the following Lemma we describe the poles of the particular meromorphic solution to problem (\ref{de.}), (\ref{icheck{v}_1^1}) constructed above.

%\newpage

%\newpage

\newpage

\subsection{Solution of difference equation, case $\Phi\neq\ds\frac{3\pi}{2}$}

We want to modify $\hat{a}_1$ into $\hat{v}_1^1$ which will satisfy all the conditions of Proposition\;\ref{pns}. To this end for $\Phi>\ds\frac{3\pi}{2}$ we add first $T_1$ to $\hat{a}_1$ removing the poles $q_1$ and $-q_1+\pi i$ (see (\ref{rcheck_1^1})), since by this proposition $\hat{v}_1^1$ must be analytic at these points. It turns out that it is possible to construct $T_1$ in such a way that it produces the pole $-p_1-\pi i$ with the desired residue, as the same proposition requires. \\
Second, we add $T_2$ producing the pole $p_1$ with the desired residue according to the same proposition.\\    
Consider 
\begin{equation}\label{T_1}
T_1(w)=\ds\frac{\pi}{2\Phi} r_1\Bigg(\coth\Big(\ds\frac{\pi(w-q_1)}{2\Phi}\Big)+\coth\Big(\ds\frac{\pi(-w+\pi i-q_1)}{2\Phi}\Big)\Bigg),
\end{equation}
where $r_1$ is given by (\ref{u2}), and $q_1$ is defined in Lemma\;\ref{Pol}.\\
It is easy to see that the function $T_1$ satisfies the following conditions: 
\begin{equation}\label{PT_1}
T_1(-w+\pi i)=T_1(w),\qquad T_1(w+2i\Phi)=T_1(w).
\end{equation}
The poles of $T_1$ belonging to $\overline{V_{-\Phi}^{\pi+\Phi}}$ are
\begin{equation}\label{pt1}
q_1,\quad -p_1-\pi i,\quad -q_1+\pi i,\quad p_1+2\pi i,
\end{equation}
(see Fig.\;\ref{Phi+1}), and
\begin{equation}\label{resw_1}
\overset{}{\underset{q_1}{res}}\;T_1=\overset{}{\underset{-p_1-\pi i}{res}}\;T_1=r_1;\qquad \overset{}{\underset{p_1+2\pi i}{res}}\;T_1=\overset{}{\underset{p_1+2\pi i}{res}}\;T_1=-r_1.
\end{equation}
Further, we define 
\begin{align}\label{T_2}
T_2(w)&:=\ds\frac{\pi}{2\Phi} r_2\Bigg(\coth\Big(\ds\frac{\pi(w-p_1)}{2\Phi}\Big)+\coth\Big(\ds\frac{\pi(-w+\pi i-p_1)}{2\Phi}\Big)\Bigg),
\end{align} 
where $r_2$ is defined by (\ref{r1r2}).
Obviously $T_2$ also satisfies (\ref{PT_1}).\\
The poles of $T_2$ in $\overline{V_{-\Phi}^{\pi+\Phi}}$ are only 
\begin{equation}\label{polT_2}
p_1,\quad -p_1+\pi i \quad {\rm and }\quad \overset{}{\underset{p_1}{res}}\;T_2=r_2,\quad \overset{}{\underset{-p_1+\pi i}{res}}\;T_2=-r_2.
\end{equation}
Finally, we define
\begin{equation}\label{hat{v}_1^1}
\hat{v}_1^1(w):=\left\{
\begin{array}{rcl}
\hat{a}_1(w)+T_1(w)+T_2(w),& & \Phi>\ds\frac{3\pi}{2}\\\\
\hat{a}_1(w)+T_2(w),& & \Phi<\ds\frac{3\pi}{2},
\end{array}
\right.
\end{equation}
where $\hat{a}_1(w)$ is given in Proposition\;\ref{P_2a_2}.

\newpage

\begin{teo}\label{Theo2} 
Let $\Phi\neq\ds\frac{3\pi}{2}$.\\ 
{\bf i)} The function $\hat{v}_1^1$ satisfies all the hypothesis of Proposition\;\ref{pns}.\\
{\bf ii)} $\hat{v}_1^1(w)\in\mathcal{H}\Big(\overline{V_{\pi}^{\frac{3\pi}{2}}}\setminus\Gamma_{\pi}\Big)$ and it has a unique pole $-p_1+\pi i$ on $\Gamma_{\pi}$  with residue 
\begin{align}\label{res}
\overset{}{\underset{-p_1+\pi i}{res}}\;\hat{v}_1^1=-r_2.
\end{align}
{\bf iii)} $\hat{v}_1^1(w)\in\mathcal{M}\Big(\overline{V_{-\frac{\pi}{2}-\Phi}^{-\Phi}}\Big)$ and it has a unique pole at $p_1-2\pi i$ for $\Phi>\ds\frac{3\pi}{2}$, 
\begin{equation}\label{ren1}
\overset{}{\underset{p_1-2\pi i}{res}}\;\hat{v}_1^1=r_2.
\end{equation}
 \end{teo}
{\bf Proof.} {\bf i)} $\hat{v}_1^1$ satisfies (\ref{de.}) and (\ref{icheck{v}_1^1}) by Proposition\;\ref{P_2a_2}, (\ref{hat{v}_1^1}) and (\ref{PT_1}).\\
Let us prove (\ref{cu1}) and (\ref{u2}).
Consider $\Phi>\ds\frac{3\pi}{2}$. By Proposition\;\ref{P_2a_2}, (\ref{resw_1}), (\ref{pt1}) and  (\ref{polT_2}), the possible poles of $\hat{v}_1^1$ in $\overline{V_{-\Phi}^{\pi+\Phi}}$ belong to  $\Big\lbrace q_1, -q_1+\pi i\Big\rbrace \cup P$, where $P$ is given by (\ref{P}).
Moreover, $\overset{}{\underset{q_1}{res}}\;\hat{v}_1^1=\overset{}{\underset{-q_1+\pi i}{res}}\;\hat{v}_1^1=0$ by (\ref{rcheck_1^1}), (\ref{resw_1}) and (\ref{polT_2}). Hence, $\hat{v}_1^1$ satisfies (\ref{cu1}). Moreover, by (\ref{resw_1}), 
$\overset{}{\underset{-p_1-\pi i}{res}}\;\hat{v}_1^1=r_1$, $\overset{}{\underset{p_1}{res}}\;\hat{v}_1^1=r_2$. Thus, (\ref{u2}) is proven for
 $\Phi>\ds\frac{3\pi}{2}$, and $\hat{v}_1^1$ satisfies all the hypothesis of Proposition\;\ref{pns} in this case (see Fig.\;\ref{phi+}).\\
Consider $\Phi<\ds\frac{3\pi}{2}$. In this case all the poles of $\hat{v}_1^1$ in $\overline{V_{-\Phi}^{\pi+\Phi}}$ belong to $P$ by Proposition\;\ref{P_2a_2}, (\ref{hat{v}_1^1}). Hence, (\ref{cu1}) holds. Moreover, from (\ref{3.4''}), (\ref{polT_2})
$
\overset{}{\underset{p_1+2\pi i}{res}}\;\hat{v}_1^1=\overset{}{\underset{p_1+2\pi i}{res}}\;\hat{a}_1=-r_1.
$
Hence,  $\overset{}{\underset{p_1-\pi i}{res}}\;\hat{v}_1^1=r_1$ by (\ref{icheck{v}_1^1}).
The equality $\overset{}{\underset{p_1}{res}}\;\hat{v}_1^1=r_2$ follows from (\ref{hat{v}_1^1}), (\ref{polT_2}) and the analyticity of $\hat{a}_1$ in $\hat{\Pi}$. Thus $\hat{v}_1^1$ satisfies (\ref{u2}) too (see Fig.\;\ref{Phi+1}).\\
{\bf ii)} By Proposition\;\ref{P_2a_2}, $a_1\in\mathcal{H}\Big(\overline{V_{-\frac{\pi}{2}-\Phi}^{\frac{\pi}{2}+\Phi}}\Big)$ which  implies $a_1\in\mathcal{H}\Big(\overline{V_{\pi}^{\frac{3\pi}{2}}}\Big)$. By (\ref{pt1}),  $T_1$ has poles in $\overline{V_{-\Phi}^{\pi+\Phi}}$ only at $q_1$, $-p_1-\pi i$, $-q_1+\pi i$, $p_1+2\pi i$. For $\Phi>\ds\frac{3\pi}{2}$ none of these poles belong to $\overline{V_{\pi}^{\frac{3\pi}{2}}}$. $T_2$ has a pole in $\overline{V_{-\Phi}^{\pi+\Phi}}$ only in $p_1$, $-p_1+\pi i$ by (\ref{T_2}). Further, $p_1\notin\overline{V_{\pi}^{\frac{3\pi}{2}}}$, $-p_1+\pi i\in\Gamma_{\pi}$ and hence ${\bf ii)}$ holds by (\ref{hat{v}_1^1}), (\ref{polT_2}).\\  
{\bf iii)} Consider $\overline{V_{-\frac{\pi}{2}-\Phi}^{-\Phi}}$. By Proposition\;\ref{P_2a_2}, $\hat{a}_1$ has poles here at $p_1-2\pi i$ for $\Phi>\ds\frac{3\pi }{2}$ and at $-p_1+\pi i-2\Phi i$ for $\Phi<\ds\frac{3\pi}{2}$ with  residues (\ref{rcheck_1^1}) and (\ref{3.4''}).\\  
From (\ref{T_1}) and  (\ref{T_2}) it follows that $T_1$ does not have poles at $\overline{V_{-\frac{\pi}{2}-\Phi}^{-\Phi}}$ and $T_2$ has a unique pole at $-p_1+\pi i-2i\Phi$ here only for $\Phi\leq\ds\frac{3\pi}{2}$ and $\overset{}{\underset{-p_1+\pi i-2i\Phi}{res}}\;T_2=-r_2$.
Hence $\hat{v}_1^1$ has a pole at $p_1-2\pi i$ for $\Phi> \ds\frac{3\pi}{2}$ and a possible pole at $-p_1+\pi i-2i\Phi$ for $\Phi< \ds\frac{3\pi}{2}$.\\
From (\ref{hat{v}_1^1}) we obtain
$$
\overset{}{\underset{p_1-2\pi i}{res}}\;\hat{v}_1^1=\overset{}{\underset{p_1-2\pi i}{res}}\;\hat{a}_1+\overset{}{\underset{p_1-2\pi i}{res}}\;T_1+\overset{}{\underset{p_1-2\pi i}{res}}\;T_2=r_2+0+0=r_2.
$$
Similarly 
$$
\overset{}{\underset{-p_1+\pi i-2i\Phi}{res}}\;\hat{v}_1^1=\overset{}{\underset{-p_1+\pi i-2i\Phi}{res}}\;\hat{a}_1+\overset{}{\underset{-p_1+\pi i-2i\Phi}{res}}\;T_1+\overset{}{\underset{-p_1+\pi i-2i\Phi}{res}}\;T_2=r_2+0-r_2=0.
$$
Therefore {\bf iii)} and  hence Theorem\;\ref{Theo2} are proven. ~~~~~$\blacksquare$

\begin{coro}\label{cv1}
For $\Phi\neq\ds\frac{3\pi}{2}$ a unique pole of $\hat{v}_1$ belonging to $\overline{V_{-\frac{\pi}{2}-\Phi}^{\frac{3\pi}{2}}}$ is $-p_1+\pi i$ and
\begin{equation}\label{reshat{v}1}
\overset{}{\underset{-p_1+\pi i}{res}}\;\hat{v}_1=2i\sin\Phi.
\end{equation}
\end{coro}
{\bf Proof.} A unique pole of $\hat{v}_1^1(w)$ in $\overline{V_{\pi}^{\frac{3\pi}{2}}}$ is only $-p_1+\pi i$ by Theorem\;\ref{Theo2} ii) with residue (\ref{res}). Hence, $\hat{v}_1$ has a unique pole at $-p_1+\pi i$ in $\overline{V_{\pi}^{\frac{3\pi}{2}}}$ and (\ref{reshat{v}1}) follows from (\ref{Si2}), (\ref{res}) and (\ref{r1r2}).
The function $\hat{v}_1$ is analytic in $\hat{V}_{\Sigma}=V_{-\Phi}^{\pi}$ by Proposition\;\ref{pns}, since $\hat{v}_1^1$ satisfies all the hypothesis of this Proposition by Theorem\;\ref{Theo2}.\\ 
It remains only to prove that $\hat{v}_1$ is analytic in $\overline{V_{-\frac{\pi}{2}-\Phi}^{-\Phi}}$. By Theorem\;\ref{Theo2} the function $\hat{v}_1^1$ has a unique pole at $p_1-2\pi i$ in $\overline{V_{-\frac{\pi}{2}-\Phi}^{-\Phi}}$ with  residue (\ref{ren1}) and the function $\hat{G}$ also has a unique pole at this point with residue $-r_2$ by (\ref{r1r2}). Hence, the function $\hat{v}_1$ is analytic in $\overline{V_{-\frac{\pi}{2}-\Phi}^{-\Phi}}$ for $\Phi\neq\ds\frac{3\pi}{2}$ and the Corollary is proven.~~~~~$\blacksquare$

\section{$h_1$-invariant solution of the difference equation in the case $\Phi=\ds\frac{3\pi}{2}$} 
\setcounter{equation}{0}

In the previous sections we have constructed a solution to problem (\ref{de.}), (\ref{icheck{v}_1^1}), satisfying all the conditions of Proposition\;\ref{pns} for $\Phi\neq\ds\frac{3\pi}{2}$.\\
It is possible to construct a solution for $\Phi=\ds\frac{3\pi}{2}$ using the same method. A slight technical inconvenience in this case arises from the fact that the function $\check{G}_2(t)$ has a pole on $\check{\beta}$. Nevertheless, one can obtain a solution with the properties indicated in Theorem\;\ref{Theo2}.\\
However, we prefer to find a solution of the problem in the case $\Phi=\ds\frac{3\pi}{2}$ by another method.
The point is that in this case it is easy to find a solution of the difference equation (\ref{de.}) in an {\it explicit form} without using the Cauchy-type integral.\\
Using the Liouville theorem  it is easy to show that this elementary solution coincides with the solution obtained by the Cauchy-type integral.\\   
In this section we give a meromorphic $h_1$-invariant solution of (\ref{de.}).

\subsection{Meromorphic solution of the difference equation for $\Phi=\ds\frac{3\pi}{2}$}

\begin{figure}[htbp]
\centering
\includegraphics[scale=0.25]{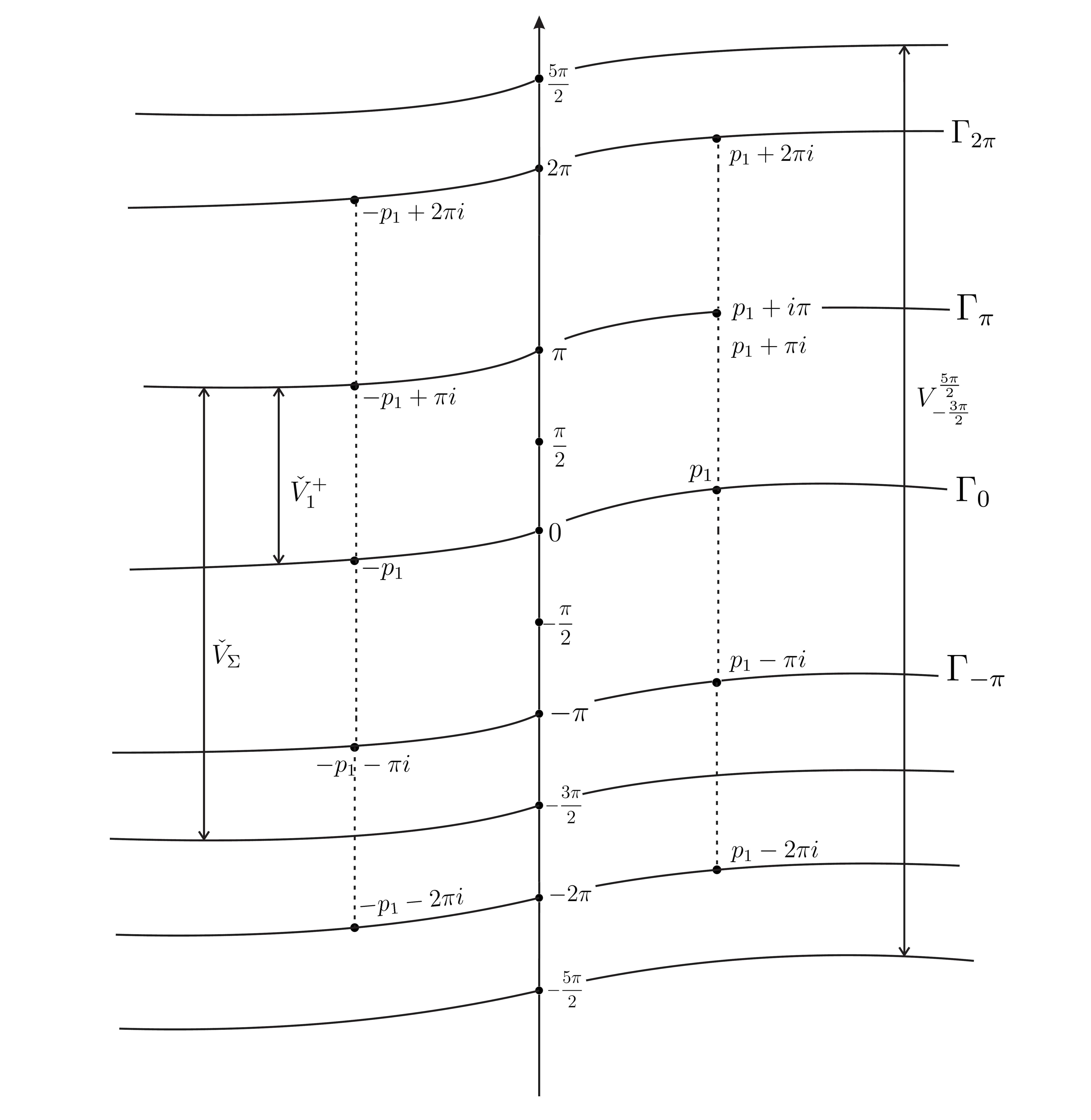}
\caption{Poles of the function $m$, $\Phi=\ds\frac{3\pi}{2}$}\label{AAJ}
\end{figure}

In this case the construction of a meromorphic solution of difference equation (\ref{de.}) is simpler than in the case $\Phi\neq\ds\frac{3\pi}{2}$ and $\hat{v}_1^1$ is expressed through elementary functions.
By (\ref{d21}), for $\Phi=\ds\frac{3\pi}{2}$, we have 
\begin{equation}\label{G_2}
\begin{array}{lll}
G_2(w) = \ds\frac{i\omega^2\sinh 2w}{\omega^2\sinh^2 w+k^2}.
\end{array}
\end{equation}
Let us solve difference equation (\ref{de.}) in this case.
First, we solve (\ref{de.}) in the class of meromorphic functions. It is easy to guess a solution, using the $3\pi i$-periodicity of $G_2$. Let us define
\begin{equation*}
m_1(w):= \frac{i w \;G_2(w)}{3\pi}.
\end{equation*}
Then, by (\ref{G_2}),  $m_1$ satisfies (\ref{de.}).
Of course, this solution is not unique. All the other solutions differ from it by a $3\pi i$-periodic function. Similarly to the case $\Phi\neq\ds\frac{3\pi}{2}$, we will modify this solution in such a way that it will satisfy all the conditions of Proposition\;\ref{pns}.\\
Function (\ref{G_2}) is not automorphic with respect to $\hat{h}_1$. Let us symmetrize it.\\
Define 
\begin{equation}\label{overline u}
m(w):=\frac{m_1(w)+m_1(-w+\pi i)}{2}.
\end{equation}
Then
\begin{equation}\label{m(w)}
m(w)=\ds\frac{\pi+2i w}{6\pi}\cdot G_2(w). 
\end{equation}
\begin{lem}\label{lu}
i) The function $m$ is an $h_1$-automorphic solution to (\ref{de.}).\\
ii) $m$ has poles in $\overline{V_{\frac{-5\pi}{2}}^{\frac{5\pi}{2}}}$ only at the points of the set 
\begin{equation}\label{p3}
P_1:=\Big\lbrace \pm p_1,  p_1\pm\pi i, -p_1\pm\pi i,  p_1\pm2\pi i, -p_1\pm2\pi i\Big\rbrace,
\end{equation}
(see Fig.\;\ref{AAJ}), and
\begin{equation}\label{re20}
\begin{array}{lll}
&&\overset{}{\underset{p_1}{res}}\;m =m_1,\quad \overset{}{\underset{p_1+\pi i}{res}}\;m =m_2,\quad \overset{}{\underset{p_1-\pi i}{res}}\;m =m_3,\quad \overset{}{\underset{p_1+2\pi i}{res}}\;m =m_4,\quad \overset{}{\underset{p_1-2\pi i}{res}}\;m =m_5,\\\\
&& \overset{}{\underset{-p_1}{res}}\;m = -m_2,\quad \overset{}{\underset{-p_1+\pi i}{res}}\;m =-m_1, \quad \overset{}{\underset{-p_1-\pi i}{res}}\;m =-m_4,\quad\overset{}{\underset{-p_1+2\pi i}{res}}\;m =-m_3,\quad \overset{}{\underset{-p_1-2\pi i}{res}}\;m=m_6,
\end{array}
\end{equation}
where
\begin{equation}\label{mk}
\begin{array}{lll}
&& m_1=-\ds\frac{p_1}{3\pi}+\ds\frac{i}{6},\quad m_2=-\ds\frac{p_1}{3\pi}-\ds\frac{i}{6},\quad m_3=-\ds\frac{p_1}{3\pi}+\ds\frac{i}{2},\quad m_4=-\ds\frac{p_1}{3\pi}-\ds\frac{i}{2},\\\\
&& m_5=-\ds\frac{p_1}{3\pi}+\ds\frac{5i}{6},\quad m_6=\ds\frac{p_1}{3\pi}+\ds\frac{5i}{6}.
\end{array}
\end{equation}
\end{lem}
{\bf Proof}. i) The assertion follows from a direct substitution of (\ref{m(w)}) into (\ref{de.}), and  (\ref{icheck{v}_1^1}) follows from (\ref{overline u}).

\medskip

 ii) The zeros of $\omega^2\sinh^{2} w+k^2$ are 
\begin{equation*}
\pm p_1+2k\pi i,\quad \pm p_1+\pi i+2k\pi i,\quad k\in\Z,
\end{equation*}
where $p_1$ is defined by (\ref{P}). Obviously, only the poles from $P_1$ belong to $\overline{V_{-\frac{5\pi}{2}}^{\frac{5\pi}{2}}}$, see Fig.\;\ref{AAJ}.\\ 
Formulas (\ref{re20}) follow from (\ref{m(w)}) and Lemma\;\ref{lh_1}.~~~~$\blacksquare$
\medskip

Now we modify the function $m(w)$ in such a way that it will satisfy the conditions (\ref{cu1}), (\ref{u2}) of Proposition\;\ref{pns}.
To this end we add to $m$ an appropriate $3\pi i$-periodic function.\\
Since for $\Phi=\ds\frac{3\pi}{2}$, $r_1=r_2=i$, conditions (\ref{cu1}), (\ref{u2})  take the form
\begin{equation*}
\hat{v}_1^1(w)\in\mathcal{H}\Big(V_{-\frac{3\pi}{2}}^{\frac{5\pi}{2}}\setminus P\Big),
\end{equation*}
where $P$ is given by (\ref{P}), and
$
\overset{}{\underset{p_1}{res}}\;\hat{v}_1^1(w)=i, \overset{}{\underset{-p_1-\pi i}{res}}\;\hat{v}_1^1(w)=i.
$

%\newpage

\vspace{0.5cm}
\subsection{Solution of the difference equation for $\Phi=\ds\frac{3\pi}{2}$}

By (\ref{p3}), the function $m$ has 8 poles in $\overline{V_{-\frac{3}{2}\pi}^{\frac{5}{2}\pi}}$ belonging to $P_1$ with residues  (\ref{re20}).
We modify $m$ so that (\ref{cu1}) and (\ref{u2}) hold.  
To this end we first add to $m$  functions $Q_1$ and $Q_3$ which ``correct'' the residues $-m_1$ at the point $p_1$ and $-m_4$ at the point $-p_1-\pi i$ by $i$. Then we add $Q_2$ which anihilates the poles $-p_1, p_1+\pi i$ and $p_1-\pi i, -p_1+2\pi i$.\\
So, consider the following functions (the periodic supplements)
\begin{align}\label{overline{T}_1}
Q_1(w)&:=\ds\frac{i-m_1}{3}\Bigg[\coth \ds\frac{w-p_1}{3}-\coth \ds\frac{w-(-p_1+\pi i)}{3}\Bigg],\\ \nonumber\\
\label{overline{T}_2}
Q_2(w)&:=-\ds\frac{m_2}{3}\Bigg[\coth \ds\frac{w-(p_1-\pi i)}{3}-\coth \ds\frac{w-(-p_1)}{3}\Bigg],
\\ \nonumber\\ \label{T_3}
Q_3(w)&:=-\ds\frac{m_3}{3}\Bigg[\coth\ds\frac{w-(p_1-\pi i)}{3}-\coth\ds\frac{w-(-p_1-\pi i)}{3}\Bigg],
\end{align}
where $m_{1,2,3}$ are given by (\ref{mk}).
Obviously, the functions $Q_{1,2,3}$ are $3\pi i$-periodic, and $\hat{h}_1$-automorphic:  
\begin{equation}\label{PT_{1,2}}
Q_{1,2,3}(w+3\pi i)=Q_{1,2,3}(w),\qquad Q_{1,2,3}\big(\hat{h}_1 w\big)=Q_{1,2,3}(w).
\end{equation}
Finally, define
\begin{equation}\label{T}
\begin{array}{lll}
Q(w)&:=&Q_1(w)+Q_2(w)+Q_3(w)
=\ds\frac{i-m_1}{3}\coth \ds\frac{w-p_1}{3}-\ds\frac{m_2+m_3}{3} \coth \ds\frac{w-(p_1-\pi i)}{3}-\\\\
&-&\ds\frac{i-m_1}{3} \coth \ds\frac{w-(-p_1+\pi i)}{3}
+\ds\frac{m_2}{3} \coth \ds\frac{w-(-p_1)}{3}+\ds\frac{m_3}{3} \coth\ds\frac{w-(-p_1-\pi i)}{3}. 
\end{array}
\end{equation}
From (\ref{overline{T}_1})-(\ref{T}), it follows directly that the  set of the poles of 
 $Q$ in $\overline{V_{-\frac{5\pi}{2}}^{\frac{5\pi}{2}}}$ is $P_1$ given by (\ref{p3}) (see Fig.\;\ref{AAJ}) and
\begin{equation}\label{PolT}
\begin{split}
&\overset{}{\underset{p_1}{res}}\;Q =i-m_1,\quad\overset{}{\underset{p_1+\pi i}{res}}\;Q=\overset{}{\underset{p_1-2\pi i}{res}}\;Q=-m_2,\quad\overset{}{\underset{p_1-\pi i}{res}}\;Q=\overset{}{\underset{p_1+2\pi i}{res}}\;Q=-m_3,\\\\
&\overset{}{\underset{-p_1-2\pi i}{res}}\;Q=\overset{}{\underset{-p_1+\pi i}{res}}\;Q=-(i-m_1), \quad\overset{}{\underset{-p_1}{res}}\;Q=m_2,\quad \overset{}{\underset{-p_1+2\pi i}{res}}\;Q=\overset{}{\underset{-p_1-\pi i}{res}}\;Q=m_3. 
\end{split}
\end{equation}

\medskip
Define
\begin{equation}\label{u_3}
\hat{v}_1^1(w):=m(w)+Q(w),
\end{equation}
where $m$ is given by (\ref{m(w)}) and $Q$ is given by (\ref{T}).
\begin{teo}\label{lu_2}
{\bf i)} For $\Phi=\ds\frac{3\pi}{2}$ the function $\hat{v}_1^1$  satisfies all the conditions of Proposition\;\ref{pns}.\\ 
{\bf ii)}  The poles of $\hat{v}_1^1$ in $\overline{V_{-\frac{5\pi}{2}}^{\frac{3\pi}{2}}}$ are
\begin{equation}\label{Pol}
-p_1-\pi i,\quad -p_1+\pi i, \quad p_1,\quad p_1-2\pi i
\end{equation} 
with the following residues
\begin{equation}\label{5.17'}
\overset{}{\underset{-p_1-\pi i}{res}}\;\hat{v}_1^1=i,\quad \overset{}{\underset{-p_1+\pi i}{res}}\;\hat{v}_1^1=-i,\quad \overset{}{\underset{p_1}{res}}\;\hat{v}_1^1=i,\quad \overset{}{\underset{p_1-2\pi i}{res}}\;\hat{v}_1^1=i.
\end{equation}
\end{teo}
\no {\bf Proof.} {\bf i)} Equations (\ref{de.}) and (\ref{icheck{v}_1^1}) follow from Lemma\;\ref{lu} and (\ref{u_3}).
From (\ref{re20}), (\ref{mk}) and (\ref{PolT}) we obtain that (\ref{u2}) holds.\\
Let us prove (\ref{cu1}). Since all the poles of $\hat{v}_1^1$ in $\overline{V_{-\frac{3\pi}{2}}^{\frac{5\pi}{2}}}$ belong to $P_1$ by Lemma\;\ref{lu}, it suffices to prove that $\hat{v}_1^1$ is analytic in $P_1\setminus P=\Big\lbrace -p_1, p_1\pm \pi i, -p_1+2\pi i\Big\rbrace$.\\ 
From (\ref{re20}), (\ref{PolT}), (\ref{icheck{v}_1^1}) and Lemma\;\ref{lh_1} we obtain 
$
\overset{}{\underset{w=-p_1}{res}}\;\hat{v}_1^1=\overset{}{\underset{p_1+\pi i}{res}}\;\hat{v}_1^1=\overset{}{\underset{-p_1+2\pi i}{res}}\;\hat{v}_1^1=\overset{}{\underset{p_1-\pi i}{res}}\;\hat{v}_1^1=0.
$
Thus (\ref{cu1}) and i) are proven.\\ 
{\bf ii)} By (\ref{p3}) and (\ref{T}) the poles of $m$ and $Q$ in $\overline{V_{-\frac{5\pi}{2}}^{\frac{3\pi}{2}}}$ are
$
\pm p_1,  p_1\pm \pi i,  \pm p_1-2\pi i,  -p_1\pm \pi i.
$
From (\ref{re20}), (\ref{mk}), (\ref{PolT}) and (\ref{u_3}) it follows that $\hat{v}_1^1(w)$ 
does not have poles at $-p_1,  p_1+\pi i, -p_1-2\pi i$ and has poles (\ref{Pol}) with residues (\ref{5.17'}).
Statement ii) is also proven.~~~~~$\blacksquare$

\medskip
Now we establish an important property of the function $\hat{v}_1$ similar to Corollary\;\ref{cv1} for the case $\Phi=\ds\frac{3\pi}{2}$.\\
We recall that this function plays the crucial role in the construction of the Sommerfeld-type representation for the solution of the main problem. This representation will be given in the following section.
\begin{coro}\label{rv1}
For $\Phi=\ds\frac{3\pi}{2}$ the function $\hat{v}_1$ given by (\ref{Si2}) has a unique pole at $-p_1+\pi i$ belonging to $\overline{V_{-\frac{5\pi}{2}}^{\frac{3\pi}{2}}}$, and 
$
\overset{}{\underset{-p_1+\pi i}{res}}\;\hat{v}_1=-2i.
$
\end{coro}
{\bf Proof.} The function $\hat{v}_1\in\mathcal{H}\big(\hat{V}_{\Sigma}\big)=\mathcal{H}\big(V_{-\frac{3\pi}{2}}^{\pi}\big)$ by Proposition\;\ref{pns} and Theorem\;\ref{lu_2}. It suffices to analyze $\overline{V_{-\frac{5\pi}{2}}^{\frac{3\pi}{2}}}\setminus\hat{V}_{\Sigma}=\overline{V_{\pi}^{\frac{3\pi}{2}}}\cup \overline{V_{-\frac{5\pi}{2}}^{-\frac{3\pi}{2}}}$.\\
First, consider $\overline{V_{\pi}^{\frac{3\pi}{2}}}$. By Theorem\;\ref{lu_2} ii), (\ref{G(w)}), and $(2\pi i)$-periodicity of $\hat{G}$, unique poles of $\hat{v}_1^1$ and $\hat{G}$ in $\overline{V_{\pi}^{\frac{3\pi}{2}}}$ are $-p_1+\pi i$  and $p_1$, and the unique pole in $\overline{V_{-\frac{5\pi}{2}}^{-\frac{3\pi}{2}}}$ of the same function is $p_1-2\pi i$ with residue (\ref{5.17'}) and (\ref{r1r2}). Hence the statement follows from (\ref{Si2}). $~~~~~\blacksquare$

\newpage

\section{Asymptotics of $\hat{v}_1$ at infinity}
\setcounter{equation}{0}

We will need to prove (\ref{SP}). For this we have to find the asymptotics of the integrand $\hat{v}_1(w)$ at infinity.

\subsection{Asymptotics of $\hat{v}_1^1$ at  infinity}

\begin{lem}\label{lashat{v}_1^1}
For any $\Phi\in(\pi,2\pi)$ the function $\hat{v}_1^1$ admits the following asymptotics:
\begin{equation}\label{as hat{v}_1^1}
\left.
\begin{array}{rcl}
\hat{v}_1^1(w)&=&\pm \ds\frac{\sin\Phi}{\Phi} \Big(w-\ds\frac{\pi i}{2}\Big)+o\Big(e^{\mp\frac{\pi}{2\Phi}w}\Big),\\\\
\ds\frac{d}{dw}\hat{v}_1^1(w)&=&\pm\ds\frac{\sin\Phi}{\Phi}+o\big(e^{\mp \frac{\pi}{2\Phi}w}\big)
\end{array}
\right| \quad \Re w\to \pm\infty.
\end{equation}

\end{lem}
{\bf Proof.}
From (\ref{T_1}) it follows that $T_1(w)$ admits the following asymptotics 
$$
T_1(w)=o\Big(e^{\mp\frac{\pi}{2\Phi}w}\Big),\quad \Re w\to\pm\infty.
$$
Similarly, $T_2(w)$ admits the same asymptotics by (\ref{T_2}) and,  by (\ref{hat{v}_1^1}), (\ref{ashat{a}_1}),  $\hat{v}_1^1$ satisfies (\ref{as hat{v}_1^1}) in the case $\Phi\neq\ds\frac{3\pi}{2}$.\\
Consider the case $\Phi=\ds\frac{3\pi}{2}$. From (\ref{u_3}), (\ref{m(w)}), (\ref{T}) it follows that the asymptotics (\ref{as hat{v}_1^1}) holds in this case too. 
Similarly,  differentiating (\ref{u_3}) we obtain (\ref{as hat{v}_1^1})~~~~~$\blacksquare$

\begin{obs}
The asymptotics of $\hat{v}_1^1$ coincide for the cases $\Phi\neq\ds\frac{3\pi}{2}$ and $\Phi=\ds\frac{3\pi}{2}$. 
\end{obs}

\subsection{Asymptotics of $\hat{v}_1(w)$}
By (\ref{Si2}),
$
\hat{v}_1(w)=\hat{v}_1^1(w)-\hat{G}(w),\quad w\in\C,
$
where $\hat{G}(w)$ is given by (\ref{G(w)}). 
Obviously, 
\begin{align*}
 \hat{G}(w)=\pm e^{\mp i\Phi}+o\big(e^{\mp\frac{\pi}{2\Phi}w}\big),\quad
 \ds\frac{d}{dw} \hat{G}(w)=o\big(e^{\mp\frac{\pi}{2\Phi}w}\big),\quad \Re w\to\pm\infty.
\end{align*}
Hence,
\begin{align}\label{as hat{v}_1}
& \hat{v}_1(w)={\rm sign}(\Re w)\cdot\ds\frac{\sin\Phi}{\Phi}\Big(w-\ds\frac{\pi i}{2}\Big) +{\rm sign} (\Re w) e^{-{\rm sign} (\Re w)wi\Phi}+o\Big(e^{{-\rm sign}(\Re w)\frac{w\pi}{2\Phi}}\Big),
\\\nonumber\\\label{as dhat{v}_1}
&\ds\frac{d}{dw}\hat{v}_1(w)={\rm sign}(\Re w)\ds\frac{\sin\Phi}{\Phi}+o\Big(e^{{-\rm sign}(\Re w)\frac{w\pi}{2\Phi}}\Big),\quad\Re w\to\pm\infty. 
\end{align}
by (\ref{as hat{v}_1^1}).~~~~~~$\blacksquare$

%\newpage
\section{Sommerfeld-type representation of solution to problem  (\ref{u_1}) }
\setcounter{equation}{0}

\begin{figure}[htbp]
\centering
\includegraphics[scale=0.17]{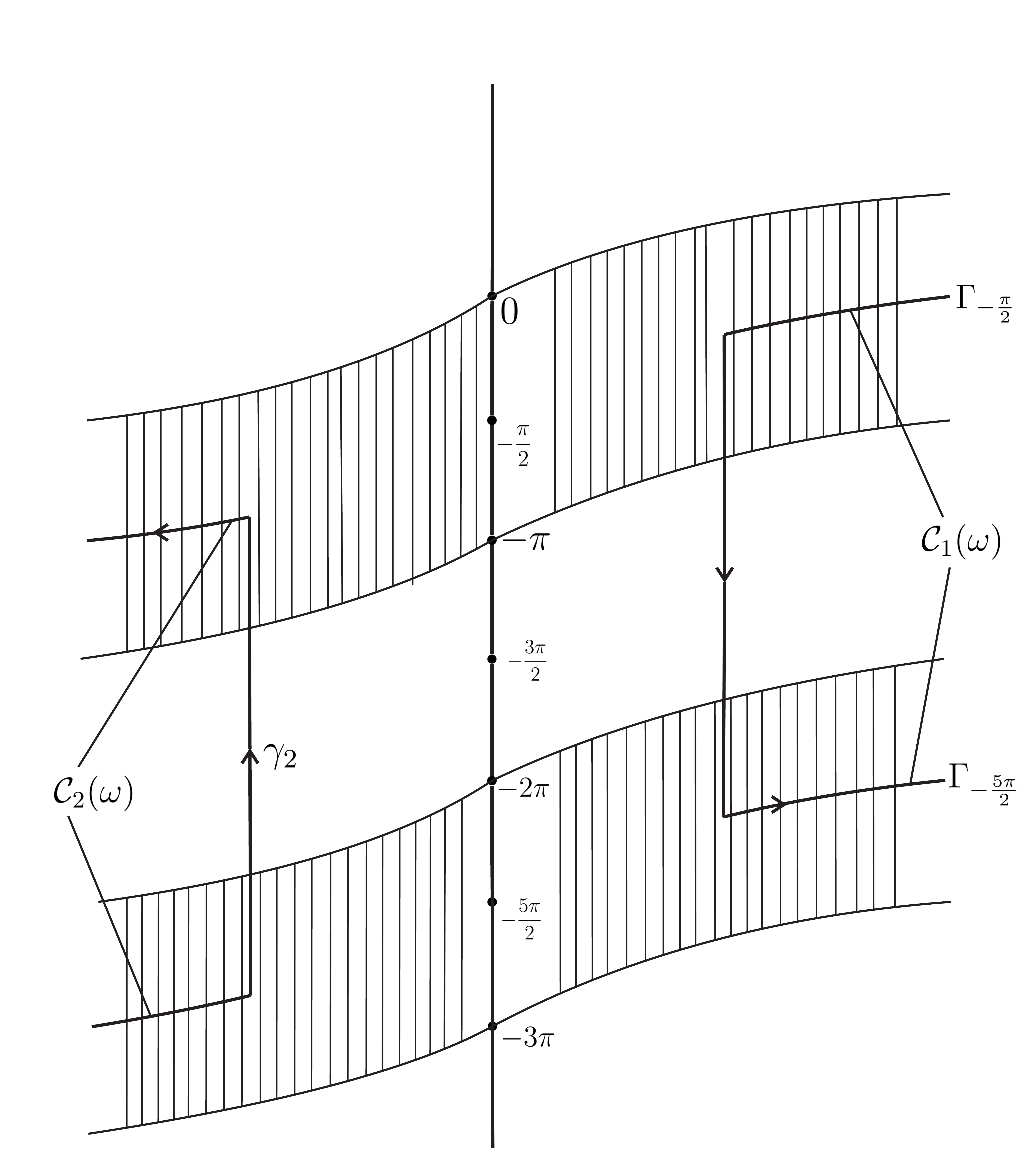}
\caption{Sommerfeld double-loop contour\;$C_{\omega}$}\label{C_1C_2}
\end{figure}

\begin{figure}[htbp]
\centering
\includegraphics[scale=0.17]{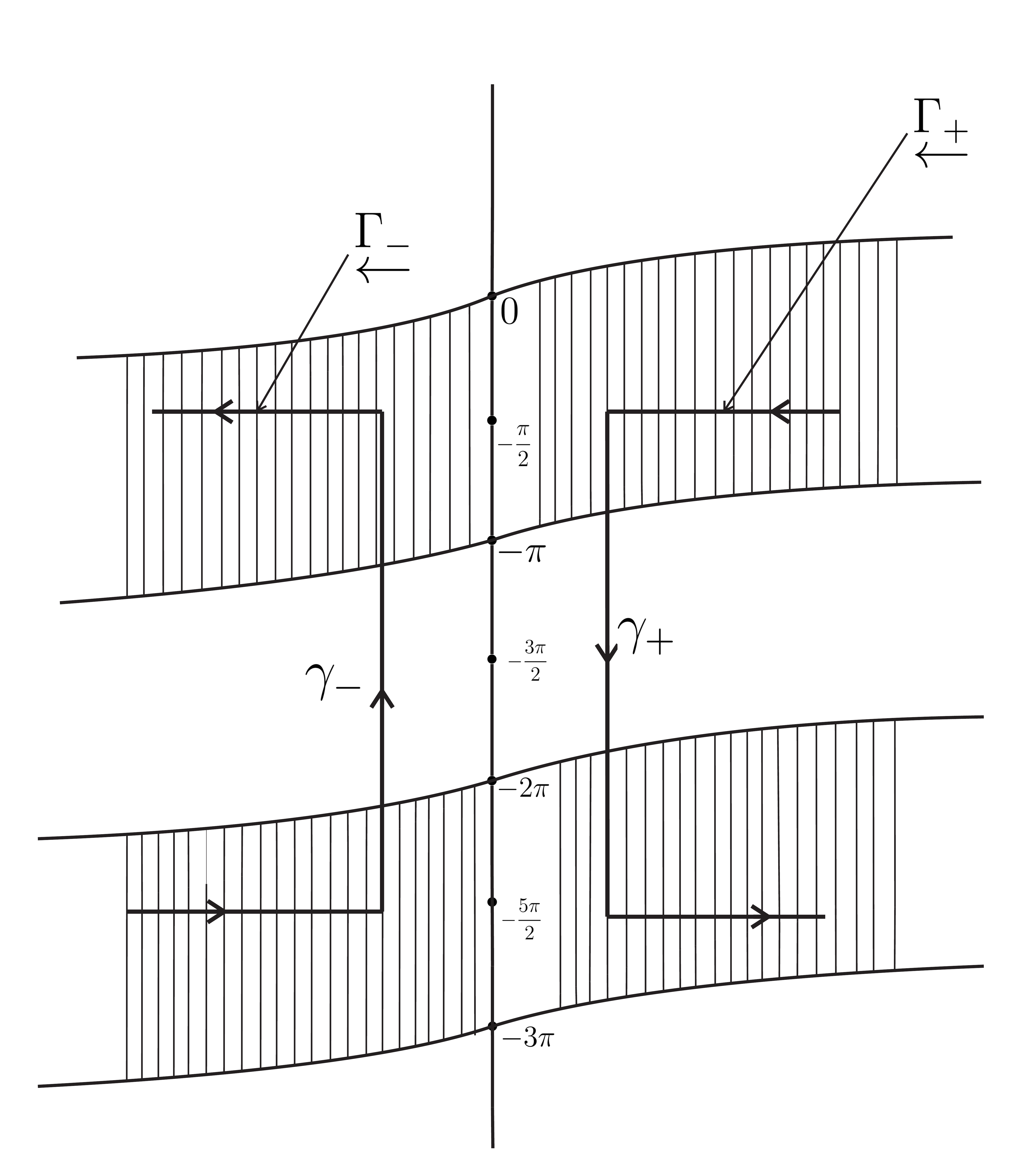}
\caption{Sommerfeld two-loop contour $\mathcal{C}=\mathcal{C}(i)$}\label{CS3}
\end{figure}

In this section we give a Sommerfeld-type representation of solution to problem (\ref{u_1}). This representation was obtained by A. Sommerfeld and it is widely used in Mathematical Diffraction Theory \cite{Somm}. This representation is an integral with a specially chosen integrand along a Sommerfeld-type contour. This contour has double-loop form as in Fig.\;\ref{CS3}.\\ 
We define first this curvilinear contour depending on $\omega\in\C^{+}$ (in contrast to the Sommerfeld contour), and then we reduce it to the rectilinear contour which coincides with the Sommerfeld contour $\mathcal{C}$ (Fig.\;\ref{CS3})).
Define $\mathcal{C}(\omega)=\mathcal{C}_1(\omega)\cup\mathcal{C}_2(\omega)$, where
\begin{equation*}
\mathcal{C}_2(\omega):=\Big\lbrace w\in\Gamma_{-\frac{5\pi}{2}}(\omega), w_1\leq -b\Big\rbrace \cup \gamma_{2}(w)\cup \Big\lbrace w\in\Gamma_{-\frac{\pi}{2}}(\omega), w_1\leq -b\Big\rbrace,
\end{equation*}  
$\gamma_{2}(\omega)$ is the segment of the line $\big\lbrace -b+iw_2, w_2\in\R\big\rbrace$ lying between $\Gamma_{-\frac{5\pi}{2}}$ and $\Gamma_{-\frac{\pi}{2}}$, $\mathcal{C}_1(\omega)=-\mathcal{C}_2(\omega)-5\pi i$ and $b\geq 2|\Re p_1|$ (see Fig.\;\ref{C_1C_2}).\\
In our case the integrand is the Sommerfeld exponential $e^{-\omega\rho\sinh w}$ multiplied by a kernel $\hat{v}_1$ which was  constructed in the previous sections.\\

\no {\bf Proof of the main Theorem\;\ref{Tu_1u_2}.}
First, we consider the Sommerfeld integral with the contour $\mathcal{C}(\omega)$. We write the integral (\ref{u_1'}) with $\mathcal{C}(\omega)$ instead of $\mathcal{C}$. We keep the notation $u_1$ for this integral because we will see later that these two integrals coincide.\\
So, let
\begin{equation}\label{2.1'}
u_1(\rho,\theta)=\ds\frac{1}{4\pi\sin\Phi}\ds\int\limits_{\mathcal{C}(\omega)} e^{-\omega\rho \sinh w}\;\hat{v}_1(w+i\theta)\;dw,
\end{equation}
where $\mathcal{C}(\omega)$ is defined above. 
Here and in the following we will use the following estimate: for $\rho>0$, $\tau\in[\tau_0,\pi-\tau_0]$ with $0<\tau_0\leq\ds\frac{\pi}{2}$, $w=w_1+iw_2\in\Gamma_0$ and $\omega\in\C^{+}$
\begin{equation}\label{eee}
\big|e^{-\omega\rho\sinh(w-i\tau)}\big|\leq e^{-C(\omega,\tau_0)\rho\cosh w_1},
\end{equation}    
where $C(\omega,\tau_0)>0$. The proof of this estimate is given in Appendix\;\ref{ic}.
Hence the integral (\ref{2.1'}) converges by  the asymptotics (\ref{as hat{v}_1}), (\ref{check{v}_1}), and (\ref{G(w)}), since $\hat{v}_1(w+i\theta)$ does not have poles on $\mathcal{C}(\omega)$ by Corollaries\;\ref{cv1} and \ref{rv1} (see Fig.\;\ref{C_1C_2}, where  the exponential decreases superexponentially in the shaded regions).\\ 
Let us prove that $u_1(\rho,\theta)$ satisfies the first equation of (\ref{u_1}). To this end we rewrite (\ref{2.1'}) as 
\begin{equation*}
u_1(\rho,\theta)=\ds\frac{1}{4\pi\sin\Phi}\ds\int\limits_{\mathcal{C}(\omega)+i\theta} e^{-\omega\rho\sinh(w-i\theta)}\; \hat{v}_1(w) dw.
\end{equation*}
Let us fix $\rho>0, \theta_0\in(\phi,2\pi)$. By the Cauchy Theorem
\begin{equation*}
u_1(\rho,\theta)=\ds\frac{1}{4\pi\sin\Phi}\ds\int\limits_{\mathcal{C}(\omega)+i\theta_0} e^{-\omega\rho\sinh(w-i\theta)}\; \hat{v}_1(w) dw
\end{equation*}
for any $\theta_0$ sufficiently close to $\theta$.\\
Now the differentiation in $(\rho,\theta)$ under the sign of the integral is possible and the first equation in (\ref{u_1}) follows from the formula $\big(\Delta+\omega^2\big) e^{-\omega\rho\sinh (w-i\theta)}=0$.\\
Finally, boundary conditions (\ref{ra2}) and (\ref{ra3}) are proved in the next section. The integral (\ref{2.1'}) is transformed into the integral (\ref{u_1'}) over the contour $\mathcal{C}=\mathcal{C}(i)$, which no longer depends on $\omega$ (see Fig.\ref{CS3}).~~~~~$\blacksquare$ 

%\newpage

\newpage

\section{Proof of the boundary conditions}
\setcounter{equation}{0}

\subsection{Decomposition of the solution into a plane wave and a wave dispersed by the vertex}

In this section we decompose the solution of problem (\ref{u_1}) given by (\ref{2.1'}) into two parts: the first part is the plane wave generated by the first boundary condition (\ref{u_1}) and the second part is the wave dispersed by the edge of the wedge.  
\newpage
\begin{figure}[htbp]
\centering
\includegraphics[scale=0.4]{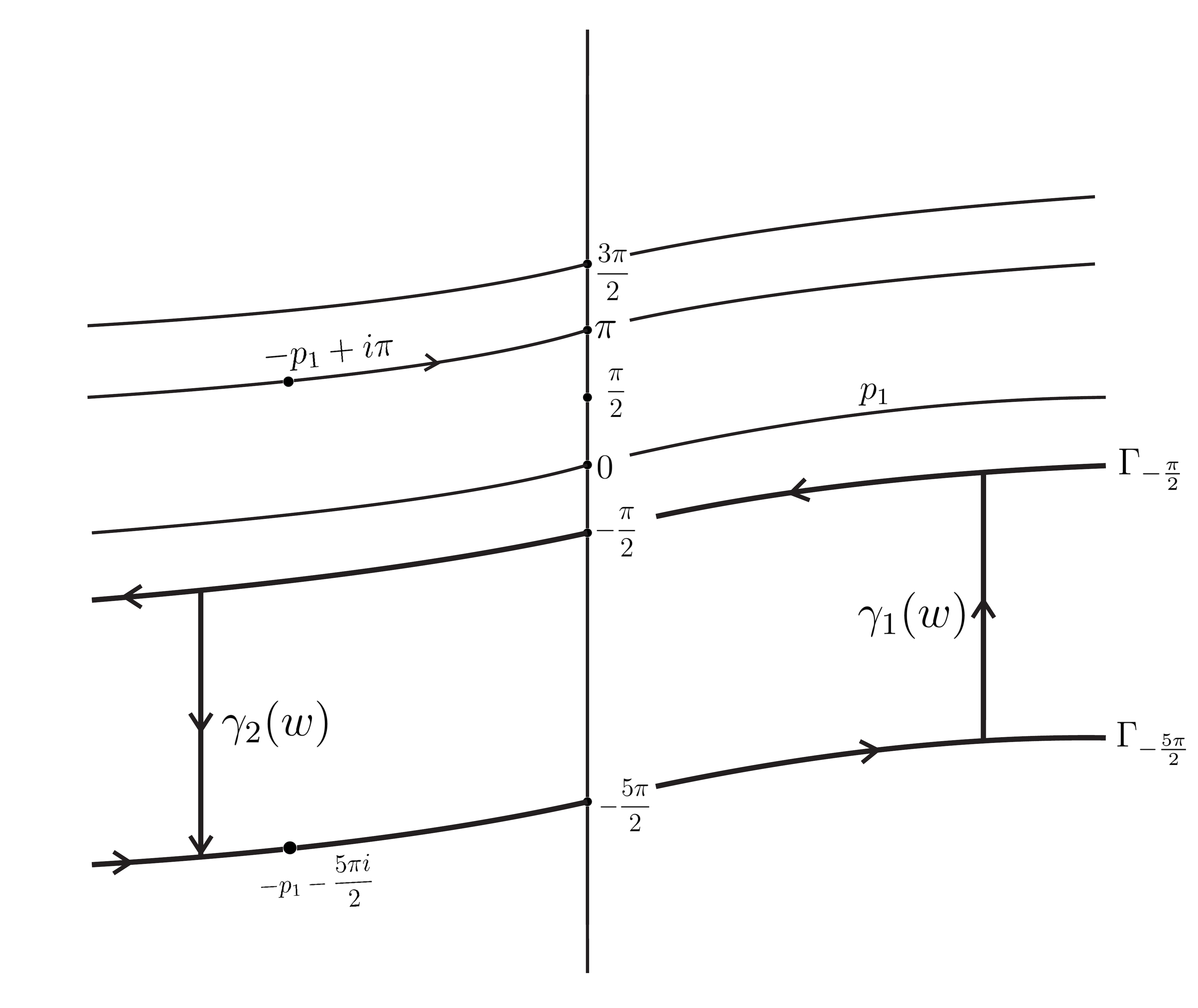}
\caption{Decomposition of the solution}\label{CC}
\end{figure}
\vspace{0.5cm}

To give this decomposition we recall that a unique pole of $\hat{v}_1(w)$ lying in $\overline{V_{-\frac{\pi}{2}-\Phi}^{\frac{3\pi}{2}}}$ is 
\begin{equation}\label{reshat{v}_1}
-p_1+\pi i\quad{\rm and}\quad \overset{}{\underset{-p_1+\pi i}{res}}\;\hat{v}_1(w)=2i\sin\Phi,
\end{equation}
as follows from Corollaries\;\ref{cv1} and \ref{rv1}.\\
Define a plane wave generated by the first boundary condition in (\ref{u_1}) ;
\begin{equation}\label{u_p}
u_{p}(\rho,\theta):=e^{-\omega\rho \sinh(p_1+i\theta)},\quad \rho>0,\quad \theta\in\R,
\end{equation}
where $p_1$, is given by (\ref{P}), and a ``diffracted'' wave 
\begin{equation*}
u_{d}(\rho,\theta):=\ds\frac{1}{4\pi\sin\Phi}\ds\int\limits_{\underrightarrow{\Gamma_{-\frac{5\pi}{2}}}\;\cup\; \underleftarrow{\Gamma_{-\frac{\pi}{2}}}} \; e^{-\omega\rho\sinh w} \hat{v}_1(w+i\theta) dw,\quad w\in(2\pi-\Phi, 2\pi].
\end{equation*}
The integrand here coincides with the integrand in (\ref{u_1'}), but the contour of integration  differs from $\mathcal{C}$ (see Fig.\;\ref{CC}).\\
It turns out that the solution $u_1$ is decomposed into the sum of $u_{p}$ and $u_{d}$ and the corresponding decomposition is more convenient for the proof of boundary conditions. 
\begin{teo}\label{td}
The solution of problem (\ref{u_1}) given by (\ref{u_1'}) admits the following representation
\begin{equation}\label{des}
u_1(\rho,\theta)=\left\{     \begin{array}{rcl}
                \ds\frac{1}{4\pi\sin\Phi}\ds\int\limits_{\underrightarrow{\Gamma_{-\frac{5\pi}{2}}}\;\cup\; \underleftarrow{\Gamma_{-\frac{\pi}{2}}}} \; e^{-\omega\rho\sinh w} \hat{v}_1(w+i\theta) dw,~~~~~~ \theta\in\Big[2\pi-\Phi,\ds\frac{3\pi}{2}\Big)\\\\ 
 \ds\frac{1}{4\pi\sin\Phi}\ds\int\limits_{\underrightarrow{\Gamma_{-\frac{5\pi}{2}}}\;\cup\; \underleftarrow{\Gamma_{-\frac{\pi}{2}}}} \; e^{-\omega\rho\sinh w} \hat{v}_1(w+i\theta) dw+u_{p}(\rho,\theta),  \theta\in\Big(\ds\frac{3\pi}{2},2\pi\Big].               
\end{array} 
\right. 
\end{equation}
\end{teo}
{\bf Proof.} By the Cauchy Theorem, $u_1$ defined by (\ref{2.1'}) admits the representation
\begin{equation}\label{8.2'}
\begin{array}{lll}
u_1(\rho,\theta)=\ds\frac{1}{4\pi \sin\Phi} \ds\int\limits_{\Gamma_{-\frac{5\pi}{2}}\cup\Gamma_{-\frac{\pi}{2}}} \;e^{-\omega\rho\sinh w} \hat{v}_1(w+i\theta) dw-\ds\frac{1}{4\pi \sin\Phi}\ds\int\limits_{\gamma(\omega)} e^{-\omega\rho\sinh w} \;\hat{v}_1(w+i\theta) dw,
\end{array}
\end{equation}
where $\gamma(\omega)$ is the contour bounded by $\gamma_{2}(\omega)$, $\gamma_1(\omega)$ and $\Gamma_{-\frac{\pi}{2}}$, $\Gamma_{-\frac{5\pi}{2}}$.\\
Let us find the poles of $\hat{v}_1(w+i\theta), w\in\Omega$ for any $\theta\in(2\pi-\Phi, 2\pi)$, where $\Omega$ is the region  bounded by $\gamma(\omega)$.\\
Let $w^{\ast}(\theta)\in\Omega$ be a pole of $\hat{v}_1(w+i\theta)$. Then $w_{p}:=w^{\ast}+i\theta$ is a pole of $\hat{v}_1(w)$ belonging to $V_{-\frac{\pi}{2}-\Phi}^{\frac{3\pi}{2}}$.
The function $\hat{v}_1(w)$ has a unique pole $w_p=-p_1+\pi i$ in $\overline{V_{-\frac{\pi}{2}-\Phi}^{\frac{3\pi}{2}}}$ with the residue (\ref{reshat{v}_1}).\\ 
Hence, $w^{\ast}(\theta)=w_{p}-i\theta=-p_1+\pi i-i\theta$. Obviously $w^{\ast}(\theta)\in\Omega$ 
only for $\theta\in\Big(\ds\frac{3\pi}{2}, 2\pi\Big]$.
Calculating the second integral in (\ref{8.2'}) with the help of residues, we obtain (\ref{des}) for $\theta\in\Big(\ds\frac{3\pi}{2}, 2\pi\Big]$.  
Therefore, (\ref{des}) holds.~~~~~$\blacksquare$

\begin{obs}
It may seem that this formula gives a discontinuous solution on the ray $\theta=\ds\frac{3\pi}{2}$, but this is not the case, since (\ref{des}) coincides with (\ref{2.1'}) by construction. Nevertheless, we also give  in Appendix\;\ref{AC} an independent proof of the continuity of $u_1(\rho,\theta)$.
\end{obs}

\newpage
\subsection{Boundary values of the solution}
We continue to prove the main theorem.
\begin{prop}
The solution $u_1(\rho,\theta)$ given by (\ref{u_1'}) is a solution to (\ref{u_1}).
\end{prop}
{\bf Proof.}
The fact that $u_1$ satisfies the Helmholtz equation in (\ref{u_1}) has been proven in Section\;6.
Let us prove the boundary conditions in (\ref{u_1}). First, we prove the first condition (\ref{u_1}) which in the polar coordinates takes the form  
$$
u_1(\rho,2\pi)=e^{-ik\rho},\quad \rho>0.
$$
Since by (\ref{u_p}), (\ref{P}),  $u_{p}(\rho,2\pi)=e^{-i k\rho}$ it suffices to prove that
$u_{d}(\rho,\theta)$ satisfies the homogeneous conditions. 
From (\ref{Si2}) we have
\begin{equation*}
\begin{array}{lll}
u_d(\rho,2\pi)
=\ds\frac{1}{4\pi\sin\Phi}\ds\int\limits_{\underrightarrow{\Gamma_{-\frac{5\pi}{2}}}\;\cup\; \underleftarrow{\Gamma_{-\frac{\pi}{2}}}} \; e^{-\omega\rho\sinh w}\Big[\hat{v}_1^1(w+2\pi i)-\hat{G}(w+2\pi i)\Big]dw,
\end{array}
\end{equation*}
since $\hat{G}$ is $2\pi i$-periodic, the integral of the second summand is equal to $0$
because $\underleftarrow{\Gamma_{-\frac{5\pi}{2}}}=-\underrightarrow{\Gamma_{\frac{\pi}{2}}}+2\pi i$.
Thus, it suffices to prove that
\begin{equation}\label{o}
\ds\int\limits_{\underrightarrow{\Gamma_{-\frac{5\pi}{2}}}\;\cup\; \underleftarrow{\Gamma_{-\frac{\pi}{2}}}} \; e^{-\omega\rho\sinh w}\hat{v}_1^1(w+2\pi i) dw=0.
\end{equation}
Making the change of the variable
$
w'=w+2\pi i,  w'\in \underrightarrow{\Gamma_{-\frac{\pi}{2}}}\;\cup\; \underleftarrow{\Gamma_{\frac{3\pi}{2}}},
$
we obtain (\ref{o}) since (after the change) the integrand is  an $h_1$-automorphic function and 
$
h_1\underrightarrow{\Gamma_{-\frac{\pi}{2}}}=\underleftarrow{\Gamma_{\frac{3\pi}{2}}}.
$\\
Let us prove the second boundary condition in (\ref{u_1}), (\ref{des}):
$
 u_1(\rho,2\pi-\Phi)=0, \rho>0.
$
From (\ref{CE}), (\ref{check v_1}) (see Remark\:\ref{r3.7}) we have
\begin{equation*}
\begin{array}{lll}
u_1(\rho,2\pi-\Phi)
=-\ds\frac{1}{4\pi\sin\Phi}\ds\int\limits_{\underrightarrow{\Gamma_{-\frac{5\pi}{2}}}\;\cup\; \underleftarrow{\Gamma_{-\frac{\pi}{2}}}} \; e^{-\omega\rho\sinh w}\; \hat{v}_2^1(w+2\pi i-i\Phi) dw.
\end{array}
\end{equation*}
Making the change of the variable $w'=w+2\pi i-i\Phi$, we obtain
$$
\\\\
 u_1(\rho,2\pi-\Phi)=-\ds\frac{1}{4\pi\sin\Phi}\ds\int\limits_{\underrightarrow{\Gamma_{-\frac{\pi}{2}-\Phi}}\;\cup\; \underleftarrow{\Gamma_{-\frac{3\pi}{2}-\Phi}}}\; e^{-\omega\rho\sinh(w-2\pi i+i\Phi)}\; \hat{v}_2^1(w) dw=0,
$$
since $\hat{v}_2^1$ is an $h_2$-automorphic function and
$
h_2\underrightarrow{\Gamma_{-\frac{\pi}{2}-\Phi}}=\underleftarrow{\Gamma_{\frac{3\pi}{2}-\Phi}}.~~~~~~\blacksquare
$ 

\section{The solution belongs to the functional class $E$ and is unique}
\setcounter{equation}{0}

\subsection{Behavior at infinity}

\begin{lem}\label{lb}
The solution $u_1$ is an $C^{\infty}$-function in $\overline{Q}\setminus \lbrace 0\rbrace$, bounded in $\overline{Q}\cap \Big\lbrace (\rho,\theta), \rho\geq\varepsilon>0\Big\rbrace$ with all its first derivatives.
\end{lem}
{\bf Proof.} This follows from the superexponential decay of the exponential $e^{-\omega\rho\sinh w}$ in the integral (\ref{2.1'}) (see shaded region in Fig.\;\ref{C_1C_2}), analyticity of $\hat{v}_1(w)$ for $|\Re w|>|\Re p_1|$ (see Corollaries\;\ref{cv1}, \ref{rv1}) and the asymptotics (\ref{as hat{v}_1}), (\ref{as dhat{v}_1}).

\subsection{Asymptotics of the solution at the origin}

We continue to prove the main theorem, namely we prove that $u_1$ given by (\ref{u_1'}) belongs to $E$. 
It remains only to prove the asymptotics (\ref{SP}).
Let us  prove the first asymptotics in (\ref{SP}). Represent the contour $\mathcal{C}$ as
\begin{equation*}
\mathcal{C}:=\big(\mathcal{C}_{+}\cup\gamma_{+}\big)\cup\big(\mathcal{C}_{-}\cup\C\big), 
\end{equation*}
where
$$
\mathcal{C}_{+}:=\underleftarrow{\Gamma_{+}}\cup \big(\underrightarrow{\Gamma_{+}}-2\pi i\big),\quad  \mathcal{C}_{-}:=\underleftarrow{\Gamma_{-}}\cup \big(\underrightarrow{\Gamma_{-}}-2\pi i\big),
$$
and the contours $\Gamma_{\pm}, \gamma_{\pm}$ are shown in Fig.\;\ref{CS3}. Note that the ``finite'' part of the integral (\ref{u_1'}), has  a ``good'' asymptotics, since
\begin{equation}\label{as f}
\ds\int\limits_{\gamma_{+}\cup\gamma_{-}} e^{-\omega\rho\sinh w}\;\hat{v}_1(w+i\theta)\;dw=C(\theta)+C_1(\theta)\rho+O(\rho),\quad\rho\to0
\end{equation}
by (\ref{as hat{v}_1}).\\
Thus it suffices only to find the principal term of the asymptotics of the ``infinite'' part of $u_1$. Since $\sinh w$ is a $2\pi i$-periodic function and is even on $\Gamma_{-}\cup\Gamma_{+}$, we have
\begin{equation*}
\begin{array}{lll}
u_{1,i}&:=&\ds\int\limits_{\mathcal{C}_{+}\cup \mathcal{C}_{-}} e^{-\omega\rho \sinh w}\;{\rm sign} (\Re w)\Bigg(\ds\frac{\sin\Phi}{\Phi}\Big(w-\ds\frac{\pi i}{2}+i\theta\Big)+ e^{-{\rm sign} (\Re w)i\Phi}\Bigg) dw=0
\end{array}
\end{equation*}
and
$$
u_1(\rho,\theta)=\ds\int\limits_{\mathcal{C}} e^{-\omega\rho\sinh w}\;\hat{v}_1(w+i\theta) dw=C(\theta)+o(1),\quad \rho\to0~~ {\rm by}~~ (\ref{as hat{v}_1}).~~~~~~\blacksquare
$$

\medskip
Let us  prove the second asymptotics in (\ref{SP}). Using  polar coordinates (\ref{PC}), we have
\begin{align*}
& \nabla u_1(\rho,\theta)=\ds\int\limits_{\mathcal{C}} \big(\partial_{y_1},\partial_{y_2}\big) e^{-\omega \rho\sinh w}\;\hat{v}_1(w+i\theta)\;dw=\ds\int\limits_{\mathcal{C}}\Big(K_1(\rho,\theta,w), K_2(\rho,\theta,w)\Big)\hat{v}_1(w)\;dw,
\end{align*}
where
$$
K_1(\rho,\theta,w)=\ds\frac{\partial}{\partial y_1}\Bigg[e^{-\omega\rho\sinh w}\;\hat{v}_{1}(w+i\theta)\Bigg]=K_{11}(\rho,\theta,w)+K_{12}(\rho,\theta,w),
$$

\begin{align}\label{K_{11} K_{12}}
 K_{11}:=e^{-\omega \rho\sinh w}\;(-\omega \sinh w) \cos\theta\;\hat{v}_1(w+i\theta),\qquad
 K_{12}:=-\ds\frac{i}{\rho}\sin\theta\;\partial_{w} \hat{v}_1(w+i\theta);
\end{align}
and
$$
K_2(\rho,\theta,w):=\partial_{y_2}\Big[e^{-\omega\rho\sinh w}\;\hat{v}_1(w+i\theta)\Big]=K_{21}(\rho,\theta,w)+K_{22}(\rho,\theta,w),
$$
where
\begin{align*}
 K_{21}:=e^{-\omega \rho\sinh w}\;(-\omega \sinh w) \sin\theta\;\hat{v}_1(w+i\theta);\qquad
 K_{22}:=\ds\frac{i\cos\theta}{\rho}\sin\theta\;\partial_{w} \hat{v}_1(w+i\theta).
\end{align*}
It suffices to find the asymptotics for $\partial_{y_1}\hat{u}_1(\rho,\theta)$, since the asymptotics of $\partial_{y_2}\hat{u}_1$ is similar.

\medskip

By (\ref{as dhat{v}_1}),
\begin{align}\label{10.4'}
\ds\int\limits_{\mathcal{C}} K_{11}(\rho,\theta,w)\;dw=\Bigg(\ds\int\limits_{\gamma_{+}\cup\gamma_{-}}+\ds\int\limits_{\mathcal{C}_{+}\cup\mathcal{C}_{-}}\Bigg) K_{11}(\rho,\theta,w)\;dw.
\end{align}
We have
\begin{equation}\label{10.4''}
\ds\int\limits_{\gamma_{+}\cup\gamma_{-}} K_{11}(\rho,\theta,w)\;dw=C(\theta)+C_1(\theta)\rho+O(\rho^2),\quad\rho\to0
\end{equation}
similarly to (\ref{as f}). We need the folowing simple statement whose proof is given in  Appendix\;\ref{lem10}. 
\begin{lem}\label{l8}   
Let $\omega\in\C^{+}$,
\begin{equation*}
A(\rho):=\ds\int\limits_{0}^{\infty} e^{i\omega\rho\cosh w}\;\cosh w\;o(e^{-\frac{\pi}{2\Phi}w}) dw.
\end{equation*}
Then
\begin{equation}\label{as A}
A(\rho)\sim C(\omega)\rho^{-1+\frac{\pi}{2\Phi}}+O(\rho^{\frac{\pi}{2\Phi}}),\quad\rho\to0.
\end{equation}
\end{lem}

\medskip
By (\ref{as hat{v}_1}), noting that $e^{-\omega\rho \sinh w}=e^{i\omega\rho\cosh w}, w\in \mathcal{C}_{+} \cup \mathcal{C}_{-}$, using Lemma\;\ref{l8} and the arguments of  the proof of (\ref{10.4'}), we obtain
\begin{equation*}
\begin{array}{lll}
\ds\int\limits_{\mathcal{C}_{+}\cup \mathcal{C}_{-}}\;K_{11}(\rho,\theta,w)\;dw
=C\rho^{-1+\frac{\pi}{2\Phi}}+C_1+O(\rho^{-\frac{\pi}{2\Phi}}),\quad\rho\to0. 
\end{array}
\end{equation*}
Hence, using  (\ref{10.4''}),
$$
\ds\int\limits_{\mathcal{C}} K_{11}(\rho,\theta,w)\;dw=C_{11}+D_{11}\rho^{-1+\frac{\pi}{2\Phi}}+O(\rho^{\frac{\pi}{2\Phi}}),\quad \rho\to0.
$$
Consider
\begin{equation}\label{10.5'}
\begin{array}{lll}
\ds\int\limits_{\mathcal{C}} K_{12}(\rho,\theta,w)\;dw &=&\Bigg(\ds\int\limits_{\gamma_{+}\cup\gamma_{-}}+\ds\int\limits_{\mathcal{C}_{+}\cup \mathcal{C}_{-}}\Bigg)\cdot K_{12}(\rho,\theta,w)\;dw.
\end{array}
\end{equation}
Similarly to (\ref{10.4''}) and using (\ref{as hat{v}_1^1}), we have
\begin{equation}\label{9.9'}
\ds\int\limits_{\gamma_{+}\cup \gamma_{-}}\;K_{12}(\rho,\theta,w)\;dw
=\ds\frac{1}{\rho} C(\theta)+C_1(\theta)+O(\rho),\quad\rho\to0.
\end{equation}
Similarly to (\ref{10.4'}), we have 
\begin{equation*}
-\ds\frac{i}{\rho}\sin\theta\ds\int\limits_{\mathcal{C}_{+}\cup \mathcal{C}_{-}} e^{-\omega \rho\sinh w}\; \Bigg[{\rm sign} (\Re w)\;\ds\frac{\sin\Phi}{\Phi}\Bigg]\;dw=0.\\
\end{equation*}
Obviously,
\begin{equation*}
\begin{array}{lll}
-\ds\frac{i}{\rho}\sin\theta\ds\int\limits_{\mathcal{C}_{+}\cup \mathcal{C}_{-}} e^{-\omega \rho\sinh w}\;o\big(e^{-{\rm sign} (\Re w)\frac{w\pi}{2\Phi}}\big) dw
=\ds\frac{1}{\rho} C_{12}+D_{12}+o(1),\quad \rho\to0.
\end{array}
\end{equation*}
Hence, substituting the asymptotics (\ref{as hat{v}_1^1}) for $\ds\frac{d}{dw} \hat{v}_1^1$ into (\ref{K_{11} K_{12}}), we obtain from (\ref{10.5'}), (\ref{9.9'}) that
\begin{equation*}
\begin{array}{lll}
\ds\int\limits_{\mathcal{C}} K_{1,2}(\rho,\theta,w)\;dw=\ds\frac{1}{\rho} C_{12}+D_{12}+o(1),\quad \rho\to0.
\end{array}
\end{equation*}
Thus,
\begin{equation*}
\ds\frac{\partial}{\partial y_1} u_1(\rho,\theta) = \ds\frac{C_{1}}{\rho}+C_{2}+o(1),\quad\rho\to0,\quad \theta\in[2\pi-\Phi,2\pi]. 
\end{equation*}
Similar asymptotics for $\partial_{y_2} u_1$ holds and  the second asymptotics of (\ref{SP}) is proven.~~~~~$\blacksquare$

\subsection{Uniqueness}

In this section we prove Statement $ii)$ of Theorem\;\ref{Tu_1u_2}.
Obviously, it suffices to prove the uniqueness of solution of problem (\ref{ra1})-(\re{ra2}) in the same space $E$.
Let $v(x)$, $\sigma(x)$ be two solutions of problem (\ref{ra1})-(\ref{ra3}) belonging to the space $E$, and $v_{l}^{\beta}(x_{l})$, $\sigma_{l}^{\beta}(x_{l})$ be their  Cauchy data ($l=1,2; \beta=0,1$). Then $\hat{v}_1^1(w), \hat{\sigma}_1^1(w)$ are $h_1$-automorphic solutions of the difference equation (\ref{de.}) and they have the same poles and residues in $V_{-\Phi}^{\Phi+\pi}$ by Proposition\;\ref{pns}. Hence, their difference $\hat{\varphi}_1^1(w):=\hat{v}_1^1(w)-\hat{\sigma}_1^1(w)$ is an analytic solution of the homogeneous equation (\ref{de.}), that is, an entire periodic function on  $\C$.\\
Moreover, since $v_{1}^{1}(x_{1})$ and $\sigma_{1}^{1}(x_{1})$ admit the same asymptotics (\ref{SP}), $\varphi_1^1(x_1)$ also admits the asymptotics (\ref{SP}). Hence its F-L transform satisfies $\tilde{\varphi}_1^1(z_1)=-\ln z+C+o(1),~z\in\C^{+},~\Re z\to+\infty$ and hence,
\begin{equation}\label{ar}
\hat{\varphi}_{1}(w) \sim C w\; {\rm sign}(w),\quad \Re w\to \pm\infty, \quad w\in \hat{V}_{1}^{+}. 
\end{equation}
Since $\hat{\varphi}_{1}(w)$ is a periodic function with period $2\Phi i$, the asymptotics (\ref{ar}) holds for $w\in\C$. This implies that 
\begin{equation}\label{ct}
\hat{\varphi}_1(w)\equiv {\rm const}.
\end{equation} 
In fact, let us apply the conformal mapping $\overline{\hat{\Pi}^{+}}\to \C^{\ast}$ given by (\ref{t(w)}).\\
It is easy to show that 
\begin{equation*}
\check{\varphi}_1(t):= \hat{\varphi}_1\big(w(t)\big)\in\mathcal{H}\Big(\C^{\ast}\setminus\lbrace 1\rbrace\Big)
\end{equation*}
and it admits the asymptotics 
$
\check{\varphi}_1(t)\sim C\log(t-1)+C, t\to1.
$
Hence it has a removable singularity at $t=1$. By the Liouville Theorem,  this implies (\ref{ct}). Therefore,   
$
\hat{v}_1^1(w)=\hat{\sigma}_1^1(w)+const
$
and hence by (\ref{al})
$
\hat{v}_2^1(w)=\hat{\sigma}_2^1(w)+const.
$
This implies that 
$
\tilde{v}_1^1(z_1)=\tilde{\sigma}_1^1(z_1)+const, \tilde{v}_2^1(z_2)=\tilde{\sigma}_2^1(z_2)+const,  z_{1,2}\in\C^{+}
$
and 
$
v_1^1(x_1)=\sigma_1^1(x_1),~~ x_1>0;   v_2^1(x_2)=\sigma_2^1(x_2), x_2>0.
$
Since $v_1^0(x_1)=\sigma_1^0(x_1)=e^{-ik_1 x_1}, x_1>0$ and $v_2^0(x_2)=\sigma_2^0(x_2)=0$ this means that the Dirichlet and Neumann data of $(v(x)-\sigma(x))$ are zeros. Hence, 
\begin{equation*}
v(x)-\sigma(x)=0
\end{equation*} 
by the uniqueness of solution of elliptic equations.~~~~~~$\blacksquare$ 
%\vspace{1cm}

\begin{figure}[htbp]
\centering
\includegraphics[scale=0.24]{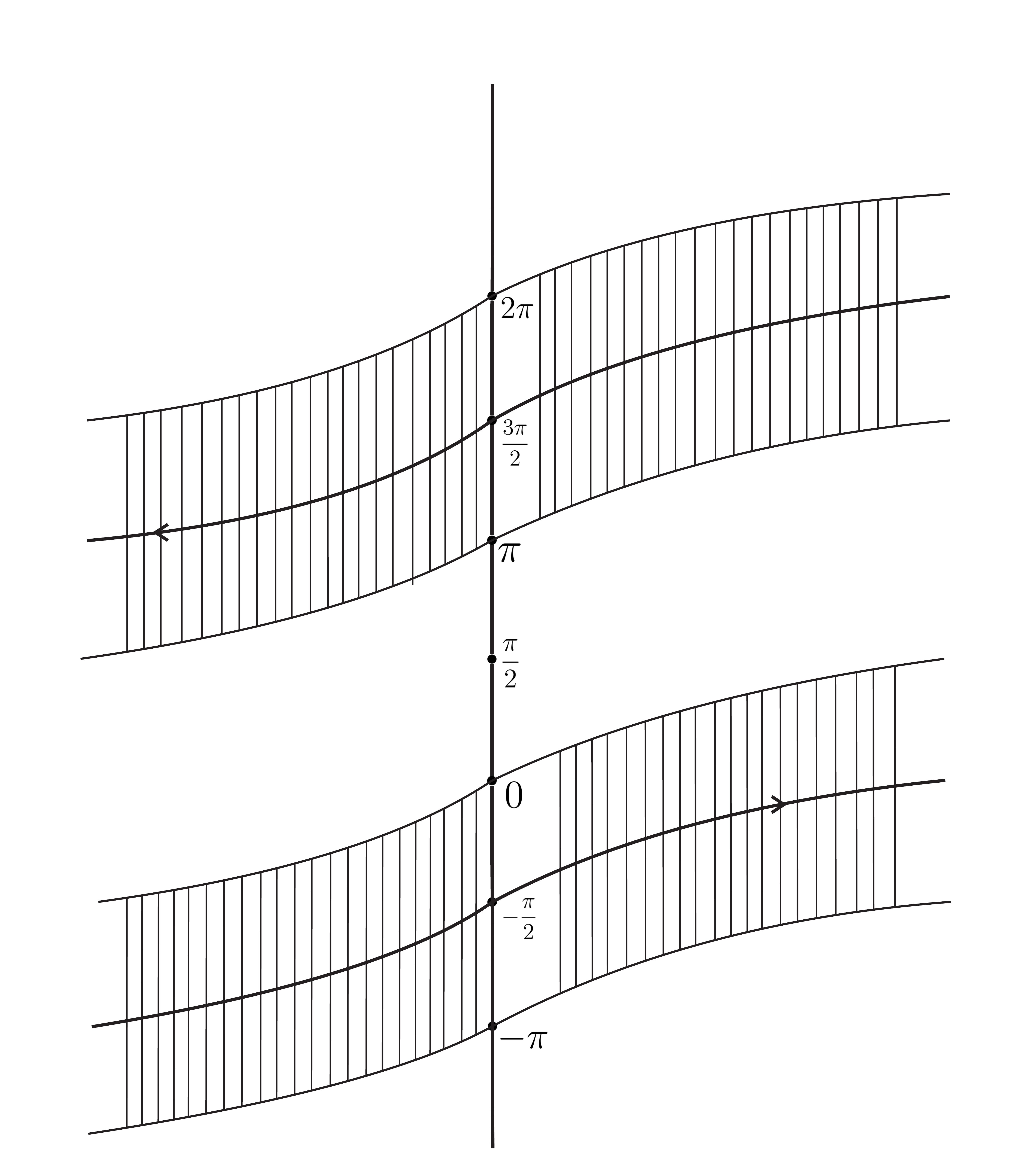}
\caption{}\label{loop33}
\end{figure}

\newpage

\section{Appendices}
\setcounter{equation}{0}

\subsection{\;\label{aal} Proof of the algebraic equation}

Since $\hat{v}_1^1(w)$ is meromorphic in $\hat{V}_{\Sigma}$ by (\ref{al}) and is  automorphic with respect to $\pi i/2$, we extend $\hat{v}_1^1$ by symmetry with respect to $\pi i/2$ to $h_1\hat{V}_{\Sigma}=\hat{V}_{0}^{\pi+\Phi}$. Namely, define 
\begin{equation*}
\overline{v}_1^1(w)=\left\{
\begin{array}{rcl}
\hat{v}_{1}^{1}(w), && w\in\hat{V}_{\Sigma}\\\\
\hat{v}_{1}^{1}(-w+\pi i),&&  w\in\hat{V}_{\pi}^{\pi+\Phi}.
\end{array}
\right.
\end{equation*}   
Obviously, $\overline{v}_1^1(w)$ is meromorphic on $\Gamma_{\pi}$ since $\hat{v}_1^1$ is meromorphic on $\Gamma_{0}$. We will still use the notation $\hat{v}_1^1$ for $\overline{v}_1^1$. Thus (\ref{icheck{v}_1^1}) holds for this extension too.\\
Since $\hat{G}(w)$ is meromorphic in $\C$, by (\ref{al}) $\hat{v}_2^1$ admits a meromorphic continuation onto $\hat{V}_{-\Phi}^{\pi+\Phi}$ which we also denote $\hat{v}_2^1$. Hence
\begin{equation}\label{3.25''}
\hat{v}_1^1(w)+\hat{v}_2^1(w)= \hat{G}(w),\quad w\in \hat{V}_{-\Phi}^{\pi+\Phi}
\end{equation} 
by the uniqueness of  analytic continuation.
Let us extend $\hat{v}_2^1(w), w\in \hat{V}_{-\pi}^{\pi+\Phi}$,  to $\hat{V}_{-\pi}^{\pi+\Phi}\cup h_2 \hat{V}_{-\pi}^{\pi+\Phi}=V_{-3\Phi}^{\pi-\Phi}$ $\Big({\rm see}~(\ref{h_1h_2})\Big)$. Similarly to $\hat{v}_1^1, \hat{v}_2^1$ is a meromorphic function in $V_{-3\Phi}^{\pi-\Phi}$. Now we extend (\ref{3.25''}) to $V_{-3\Phi}^{\pi-\Phi}$ and obtain
$
\hat{v}_1^1(w)+\hat{v}_2^1(w)=\hat{G}(w), w\in V_{-3\Phi}^{\pi-\Phi}.
$
Hence
$
\hat{v}_{l}^{1}\in \mathcal{H}(\hat{V}_{l}^{+})\cap \mathcal{M}\Big(V_{-3\Phi}^{\pi-\Phi}\Big).
$
\\ 
Continuing the process of extension of $\hat{v}_{1}^{1}$ and $\hat{v}_{2}^{1}$ by symmetries (\ref{h_1h_2}), (\ref{h_1}), we can extend equation (\ref{3.25''}) to $\C$ and obtain (\ref{3.25})-(\ref{icheck{v}_1^1}). ~~$\blacksquare$

\subsection{\label{ade}Proof of Lemma\;\ref{Vin}}

Let us apply the automorphism $h_2$ to equation (\ref{3.25}). Since by (\ref{check{v}_1}), (\ref{h_2}), 
$$
\hat{G}(h_2 w)=\ds\frac{i\omega\sinh(w+3i\Phi)}{i\omega\sinh(w+2i\Phi)+k},
$$
(\ref{al}) gives  
\begin{equation}\label{3.25'}
\hat{v}_1^1(h_2 w)+\hat{v}_2^1(h_2 w)=\hat{G}(h_2 w)=\ds\frac{i\omega\sinh(w+3i\Phi)}{i\omega\sinh(w+2i\Phi)+k},\quad w\in\C.
\end{equation}
Hence, subtracting  equation (\ref{3.25'}) from equation (\ref{al}), we obtain
$
\hat{v}_1^1(w)-\hat{v}_1^1(h_2 w)=\hat{G}_2(w),
$
where $\hat{G}_2$ is given by (\ref{d21}).\\ 
Using (\ref{h_1}), we can represent $\hat{v}_1^1\Big(h_2 w\Big)$ as a function with shifted argument. Applying (\ref{h_1h_2}) we have
$$
\hat{v}_1^1(h_2 w)=\hat{v}_1^1(h_1 h_2 w)=\hat{v}_1^1(h w),
$$
where $h(w)$ is a shift since
$
h w=h_1 h_2 w=(-h_2 w)+\pi i=(w-\pi i+2i\Phi)+\pi i=w+2i\Phi.
$
Hence, by (\ref{d21}),
$
\hat{v}_1^1(w)-\hat{v}_1^1(w+2i\Phi)=\hat{G}_2(w).~~~~~\blacksquare
$

\subsection{\label{ic}Proof of the estimate (\ref{eee})}

For $w\in\Gamma_0$, $\omega=\omega_1+i\omega_2, \omega_2>0$ consider
$
e^{-\omega\rho\sinh(w-i\tau)},\quad w\in\Gamma_0,\quad 0<\tau_0\leq\tau<\pi-\tau_0.
$
We have
%\begin{equation*}
%\begin{array}{lll}
$e^{-\omega\rho\sinh w\cos\tau}=e^{-i\omega\sinh w (-i\rho\cos\tau)}=e^{z_1(w)\cdot(-i\rho\cos\tau)}$,
%\end{array}
%\end{equation*}
where $z_1(w)$ is given by (\ref{z_{1,2}(om)}). We have for $w\in\Gamma_0$ that $z_1(w)\in\R$ by (\ref{Gamma_0'}). 
Hence $\big|e^{-\omega\rho\sinh w\cos\tau}\big|=1$ 
and,  since
$
-\omega\rho\sinh(w-i\tau)=-\omega\rho\Big[\sinh w\cos\tau-i\cosh w\sin\tau\Big],
$
we obtain
$
\big|e^{-\omega\rho\sinh(w-i\tau)}\big|=\big|e^{i\omega\rho\cosh w\sin\tau}\big|
= \big|e^{-\omega_2\rho\sin\tau\cosh w_1\cos w_2}\big| \cdot\big|e^{-\omega_1\rho\sin\tau\sinh w_1\sin w_2}\big|.
$\\
Since for $w\in\Gamma_0$, $w_2=\arctan \ds\frac{\omega_1}{\omega_2}\tanh w_1$, (see (\ref{Gamma_0'})), $\cos w_2\geq C(\omega)>0$, we have for $\tau\in[\tau_0,\pi-\tau_0]$
$$
\big|e^{-\omega_2\rho\sin\tau\cosh w_1\cos w_2}\big|\leq e^{-C(\omega,\tau_0)\rho\cosh w_1},\quad C(\omega,\tau_0)>0,\quad w_2\in\Gamma_0, \quad\omega_2>0.
$$
Moreover, for $w\in\Gamma_0$
$
\big|e^{-\omega_1\rho\sin\tau\cosh w_1\sin w_2}\big|=\Big|e^{-\omega_1^2\rho\sin\tau\;\frac{\sinh^2 w_1}{\sqrt{\omega_2^2\cosh^2 w_1+\omega_1^2\sinh^2w_2}}}\Big|\leq 1.
$
Hence, (\ref{eee}) follows.~~~~~~$\blacksquare$

\subsection{\label{AC} Analysis of the solution near the ray $\theta=\ds\frac{3\pi}{2}$}

We prove that $u_1\Big(\rho, \ds\frac{3\pi}{2}+\delta\Big)$ is continuous in $\delta$ for small  $\delta$. By (\ref{des}) it suffices to prove that 
\begin{equation}\label{I}
I(\rho,\delta):=\ds\frac{1}{4\pi\sin\Phi}\ds\int\limits_{\Gamma_{-\frac{\pi}{2}}} e^{-\omega\rho\sinh w} \hat{v}_1\Big(w+\ds\frac{3\pi i}{2}+i\delta\Big) dw
\end{equation}
satisfies 
\begin{equation}\label{SI}
I(\rho,0-)=I(\rho,0+)+u_{p}\Big(\rho,\ds\frac{3\pi}{2}\Big)
\end{equation}
since $\hat{v}_1$ does not have poles on $\Gamma_{-\pi}$ by (\ref{reshat{v}_1}).\\
Making the change of the variable $w'\to w+\ds\frac{3\pi i}{2}$ in (\ref{I}), using (\ref{reshat{v}_1}), and the Sokhotski-Plemelj Theorem, we obtain
\begin{equation*}
I(\rho,0+)=\ds\lim_{\delta\to 0}\;\ds\frac{1}{4\pi\sin\Phi}\ds\int\limits_{\underleftarrow{\Gamma_{\pi}}} e^{-\omega\rho\sinh (w-\frac{3\pi i}{2})}\;\ds\frac{\hat{v}_1(w+i\delta)\big(w-(-p_1+\pi i)+i\delta\big)}{w-(-p_1+\pi i)+i\delta} dw=-\ds\frac{1}{2} e^{i\omega\rho\cosh p_1}+V.P.,
\end{equation*}
where $V.P.$ here and in the following denotes the principal value of the integral (\ref{I}).\\
Similarly, 
$$
I(\rho, 0-)=\ds\frac{1}{2} e^{i\omega\rho\cosh p_1}+V.P.
$$
Hence (\ref{SI}) follows, since $u_{p}\Big(\rho,\ds\frac{3\pi}{2}\Big)=e^{i\omega\rho \cosh p_1}$
by (\ref{u_p}).~~~~~~$\blacksquare$

\subsection{\label{AC1} Proof of asymptotics (\ref{ashat{a}_1})}

 From (\ref{hat{a}_1(w)}), (\ref{as hat{a}_1}) it follows that
\begin{equation}\label{ashat{a}1}
\hat{a}_1(w)=-\ds\frac{\check{G}_2(1)}{2\pi i}\ln\ds\frac{1}{\coth^2\Big(\ds\frac{\pi}{2\Phi}\Big(w-\ds\frac{\pi i}{2}\Big)\Big)-1}+C+O\big(e^{\mp\frac{\pi}{\Phi}w}\big)+C_1,\quad \Re w\to\pm\infty
\end{equation}
since
\begin{equation*}
O(t(w)-1)=O\Big(e^{\mp\frac{\pi}{\Phi}w}\Big),\quad \Re w\to\pm\infty.
\end{equation*}
Let us to prove (\ref{ashat{a}_1}).
We have for $m=\ds\frac{\pi}{2\Phi}\Big(w-\ds\frac{\pi i}{2}\Big)$
\begin{equation*}
\begin{array}{lll}
\ln\ds\frac{1}{\coth^2\ds\frac{\pi}{2\Phi}\Big(w-\ds\frac{\pi i}{2}\Big)-1}
=\pm2m-\ln 4+o\big(e^{\mp m}\big),\quad \Re m\to\pm\infty;
\end{array}
\end{equation*}
\begin{equation*}
\begin{array}{lll}
-\ds\frac{\check{G}_2(1)}{2\pi i}\ln\Bigg(\ds\frac{1}{\coth^2\ds\frac{\pi}{2\Phi}\Big(w-\ds\frac{\pi i}{2}\Big)}\Bigg)&=&-\ds\frac{\check{G}_2(1)}{2\pi i}\Bigg(\pm\ds\frac{\pi}{\Phi}\Big(w-\ds\frac{\pi i}{2}\Big)-\ln 4\Bigg)+o\big(e^{\mp \frac{\pi}{2\Phi}w}\big)\\\\
&=& \pm \ds\frac{\sin\Phi}{\Phi}\Big(w-\ds\frac{\pi i}{2}\Big)-\ds\frac{\sin\Phi}{\pi} \ln 4+o\big(e^{\mp \frac{\pi}{2\Phi}w}\big).
\end{array}
\end{equation*}
Hence, using (\ref{hat{a}_1(w)}), (\ref{defC_1}), (\ref{ashat{a}1}) and (\ref{asG2}), we obtain (\ref{ashat{a}_1}).\\
Let us prove (\ref{4.27'}). From (\ref{hat{a}_1(w)}) and (\ref{t(w)}) we have 
\begin{equation*}
\begin{array}{lll}
\ds\frac{d}{dw}\hat{a}_1(w)=
 -\ds\frac{\pi}{\Phi}\;\ds\frac{\cosh\Big(\ds\frac{\pi}{2\Phi}\big(w-\ds\frac{\pi i}{2}\big)\Big)}{\sinh^3\Big(\ds\frac{\pi}{2\Phi}\big(w-\ds\frac{\pi i}{2}\big)\Big)}\;\ds\frac{d}{dt}\check{a}_1(t(w)).
\end{array}
\end{equation*}
Using (\ref{4.22'}) we obtain (\ref{ashat{a}1}).~~~~~~$\blacksquare$

\subsection{\label{lem10} Proof of Lemma\;\ref{l8}}

{\bf Proof.} Obviously, the principal term admits the asymptotics
\begin{align*}
 A(\rho)\sim \ds\int\limits_{0}^{\infty} e^{-\omega_2\rho e^{w}}\;e^{w} e^{-\frac{\pi}{2\Phi}w}\; dw,\quad \rho\to0.
\end{align*}
Making the change of the variable $\rho e^{w}=s$, we obtain
\begin{align*}
 A(\rho)\sim \rho^{-1+\frac{\pi}{2\Phi}} \ds\int\limits_{\rho}^{\infty} e^{-\omega_2 s}\;s^{(1-\frac{\pi}{2\Phi})}\;\varphi(s) ds,\quad \varphi(s)\to0,\quad s\to\infty.
\end{align*}
Hence (\ref{as A}) follows.~~~~~~$\blacksquare$

\section{Conclusion}

As is known, an angle is  one of the few regions  where  the boundary value problems for the Helmholtz equation admit an explicit solution. As far as we know, this has always been done for decreasing boundary data, where the operator methods are normally used with the exception of a very specific boundary value problem associated with the  plane incident wave (Sommerfeld's diffraction problem   \cite{Somm}). In the presented  work, we solve the Dirichlet boundary  problem not related to the incidence of a plane wave and  we  obtain an explicit  solution in the form of the Sommerfeld  integral.
The proposed method is suitable for the Neumann (NN) and Dirichlet-Neumann (DN) boundary conditions, and for angles less then $\pi$.
We hope that the method is suitable for solving such problems with a real wave number in the Helmholtz operator and also  for nonstationary  problems.

\end{document}